\documentclass[11pt]{article}
\usepackage[T1]{fontenc}
\usepackage{lmodern}
\usepackage{amsmath,amssymb,mathtools}
\usepackage[letterpaper,text={6.5in,9in}]{geometry}
\usepackage{microtype}
\usepackage[authoryear,round]{natbib}
\usepackage{paper-defs}
\usepackage{paper-style}
\title{Exact bridges for sticky reflected Brownian motion}
\author{Christopher Mount, Ibrahim Fatkullin\\{\em University of Arizona}}
\date{\today}

\begin{document}
\maketitle

\begin{abstract}
We construct exact samplers for sticky reflected Brownian motion and its local time at zero. 
Given the position and local time at both ends of an interval, we sample the resulting bridge at intermediate times in any order. 
The gain in local time between the endpoints determines the occupation time at zero. 
When the gain is positive, we draw the first and last zeros of the bridge.
The position follows a Bessel(3) bridge on each outer interval of positive length, and the entire gain in local time occurs between the two zeros.
To sample inside the middle bridge, we draw the local time at a query time by inverting an explicit distribution function.
We then sample the position by rejection with an acceptance probability bounded below uniformly in the parameters.
The construction also gives an exact forward sampler. Both samplers avoid time discretization, and we count inversion of an explicit continuous distribution function at a uniform draw as one exact operation.
Finally, we prove that the sampled values and recorded zeros determine the conditional law of the remaining path.
We can therefore refine the same path in any order.
\end{abstract}

\section{Introduction}

A Brownian particle at a sticky wall can occupy a spot on the wall for a positive amount of time even though it 
never occupies that spot for any interval of time. 
The times at which the particle is at the wall form a Cantor-like set whose positive Lebesgue measure accounts 
for the stickiness, though observations on a time grid can miss much of the time spent there.
One physical model of this behavior places the particle in a narrow, deep attractive well at the wall.
As the well narrows and deepens with its integrated Boltzmann weight tending to a finite positive limit, 
matched asymptotic expansions give the law of sticky reflected Brownian motion to leading 
order \citep{bou-rabeeStickyBrownianMotion2020}.
This limit is why sticky diffusions model particles that bind to a surface or to one another
\citep{bou-rabeeStickyBrownianMotion2020,holmes-cerfonSimulatingStickyParticles2020,bressloffCloseEncountersSticky2023}. 
Sticky reflected Brownian motion is one such diffusion, defined on $[0,\infty)$, that moves like Brownian motion away from $0$.
Feller identified its behavior at $0$ among the possible boundary conditions for a diffusion on the half line
\citep{fellerParabolicDifferentialEquations1952}, and \citet{itoBrownianMotionsHalfline1963} constructed the process from reflected Brownian motion by a random time change.
This paper will follow their construction and write $\Y_\iv$ for the resulting process, $\locT_\iv$ for its local time at $0$ and $\stick>0$ for the stickiness.
In our normalization the occupation time of $0$ up to time $s$ is $\stick\locT_s$.
(In the convention of Revuz and Yor, $L_s^0(\Y)=2\locT_s$, as in Section~\ref{sec:time-change}.)

Several applications measure this occupation time directly. In encounter-based models of surface reactions, a particle is absorbed once its occupation time at the boundary crosses a random threshold \citep{bressloffCloseEncountersSticky2023}.
In sticky piecewise deterministic samplers for sparse Bayesian models, the fraction of time a coordinate spends at zero estimates the conditional probability, given the data, that the coordinate is zero \citep{bierkensStickyPDMPSamplers2023}.
Sticky diffusions also model prices that stick at a level \citep{jiangExplicitResultsOne2019,anagnostakisPricingHedgingSticky2025}.
Estimators of the stickiness use the occupation time and the local time of an observed path \citep{anagnostakisFunctionalConvergenceLocal2023,anagnostakisNumberCrossingsBouncings2025}.
The cited studies test these estimators on grid approximations. Markov chains on a spatial grid approximate the process
\citep{bou-rabeeStickyBrownianMotion2020,meierMarkovChainApproximation2021,anagnostakisGeneralDiffusionProcesses2023}.
Weak schemes approximate stochastic differential equations with sticky boundary conditions \citep{sharmaWeakApproximationStochastic2025}, and a modified Euler--Maruyama scheme discretizes the process in time \citep{jiangModifiedEulerMaruyamaMethod2026}.
Their errors vanish as the grid or time step is refined, but these schemes do not give the exact conditional law of a new value between two observations.
\citet{shardlowStickyCIRProcess2026} samples a sticky Cox--Ingersoll--Ross process exactly at independent exponential times, which suffices for its invariant measure.
The exact algorithms for diffusions \citep{beskosExactSimulationDiffusions2005} have been extended to equations with local time
\citep{etoreExactSimulationOnedimensional2013a}, but, to our knowledge, not to a sticky boundary.
Our samplers give four quantities used in these applications. At a fixed time $s$, $\stick\locT_s$ gives
the occupation time of $0$, and $\id_{\{0\}}(\Y_s)$ records whether the path is at the boundary.
Between two fixed observation times, an increase in $\locT$ identifies a zero in the open interval almost surely.
At a positive bridge observation, the recorded zeros identify the excursion endpoints whenever they lie within the observed interval.

Our precise contributions are the following. First, Theorem~\ref{thm:forward} gives an exact forward step for the pair
$(\Y,\locT)$. It uses one normal variable, at most one uniform and at most one exponential, and no rejection. Second,
Theorem~\ref{thm:bridge} gives an exact sampler for the bridge of $(\Y,\locT)$ conditioned on its position and local
time at both ends. Third, Theorem~\ref{thm:refine} shows that a skeleton of the sampled values and zeros permits further
queries in any order, even at times computed from the skeleton, and that the result has the law of the same queries
applied to a path of the process.
When only the endpoint positions are prescribed, one first draws the terminal local time given those positions and then applies the joint-bridge sampler. The companion note \emph{Exact bridges for sticky Brownian motion with constant drift} \citep{mountFatkullinDriftNote}
gives the law and sampler for this preliminary draw and extends the construction to two-sided sticky Brownian motion
and to constant drift. Between the two zeros of a bridge, we sample the local time at a
fixed time by inverting an explicit distribution function at a single uniform. Once one has sampled the local time, 
the position has an atom at $0$ and a density on $(0,\infty)$, and a
single rejection loop samples it with an acceptance probability bounded below uniformly in the parameters
(Proposition~\ref{prop:cost}). Throughout, we count the inversion of an explicit continuous distribution function at a 
uniform as one exact draw. 

Our work rests on two classical descriptions of a Brownian path, one in terms of its running maximum and 
the other in terms of its excursions.
In a more recent work \citet{casterasTrivariateDistributionSticky2026} give the full joint law of the position, 
the local time at zero, and the occupation time of the positive half-line for sticky Brownian motion.
\citet[Thm.~3.3]{touhamiSkewStickyBrownian2021} had given only its density away from $0$. For Brownian motion, this joint
law goes back to \citet{karatzasTrivariateDensityBrownian1984}, and \citet{salminenJointDistributionBrownian1987}
explained why it matches the joint law of the maximum, its location and the endpoint. Our construction uses the same
exchange of the local time for the running maximum, and our sampler factors this law, drawing the local time first, from
a single passage time, and then the position conditionally on it. The decomposition of the bridge at its first and
last zeros is, in operational time, the decomposition of a Brownian path at its maximum
\citep{williamsPathDecompositionContinuity1974,pitmanDecompositionMaximumExcursions1996}.
Conditional on its duration, the excursion between two consecutive zeros is a Bessel(3) bridge.
The Bessel(3) process is the Doob $\mathfrak h$-transform of Brownian motion killed at $0$ with $\mathfrak h(x)=x$
\citep{doobConditionalBrownianMotion1957a,williamsPathDecompositionContinuity1974,pitmanOnedimensionalBrownianMotion1975}.

The paper is
organized as follows. In Section~\ref{sec:srbm_and_ell}, we construct $(\Y,\locT)$ from a
Brownian motion and its running maximum, and we collect the laws of passage times and Bessel(3) bridges.
Sections~\ref{sec:forward} and~\ref{sec:bridge} give the forward step and the bridge. Section~\ref{sec:refine} proves
that refinement in any order is exact. Section~\ref{sec:numerics} describes the implementation and compares the local time it draws inside a bridge with its explicit distribution function.
The appendices give the proofs.

\section{Sticky reflected Brownian motion and its local time}\label{sec:srbm_and_ell}

\subsection{Time change construction and local time normalization}\label{sec:time-change}

\begin{figure}[htb]
	\centering
	\includegraphics[width=\textwidth]{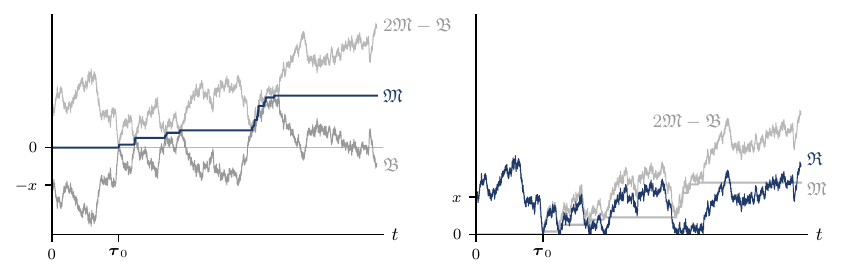}
	\caption[Reflection]{On the left, $\BM$ starts at $-x$ and $\BMM$ is its running maximum from the first passage
	$\hT_0$ to $0$ onward. Reflecting $\BM$ in $\BMM$ gives $2\BMM-\BM$. On the right, $\BMM$ is subtracted from
	$2\BMM-\BM$, which leaves the reflected Brownian motion $\BMr=\BMM-\BM$, and $\BMM$ is its local time at $0$. 
	The path is a Gaussian random walk with $10^5$ steps and $x=0.7$ (seed 8).}
	\label{fig:reflection}
\end{figure}

Let $\Omega$ be the space of continuous functions from $[0,\infty)$ to $\mathbb R$ with its Borel $\sigma$-field,
and let $\BM_\iv$ be the coordinate process. For $x\geq0$, we write $\Probx{x}$ for the Wiener measure under which
$\BM$ is a Brownian motion with $\BM_0=-x$ \citep{wiener1923differential}, and $\E_x$ for the expectation under it;
see \citet[Sec.~2.4]{karatzas1991brownian} for the construction of these measures on $\Omega$. We write $\filt_t$
for the $\sigma$-field generated by $(\BM_v)_{v\leq t}$, without completion. The proofs use only the strong Markov
property of $\BM$ at the stopping times of this filtration
(\citealp{huntTheoremsConcerningBrownian1956};\ \citealp[Thm.~2.6.15]{karatzas1991brownian}), and they hold without
change under the usual conditions. The samplers draw random variables of three kinds, namely $\uni$ uniform on
$(0,1)$, $\gauss$ standard normal and $\erv$ exponential with mean $1$, and each draw in an algorithm is a fresh
independent copy. We write $\Prob$ for the law of these draws and of what the samplers build from them. Fix
$x\geq0$. The first
passage time of $\BM$ to a level $a\geq0$ is $\hT_a\ass\inf\{t:\BM_t=a\}$. We write $\BMM_\iv$ for its running
maximum after the first passage to $0$, that is,
\[
	\BMM_t~\ass~\begin{cases}0,& t<\hT_0,\\ \displaystyle\max_{\hT_0\leq u\leq t}\BM_u,& t\geq\hT_0,\end{cases}
\]
and we set $\BMr_\iv\ass\BMM_\iv-\BM_\iv$. Since $\BMM$ is a running maximum, it first reaches a level $a>0$ at the
passage time $\hT_a$, and the samplers use this fact throughout. Then $\BMr_\iv$ is reflected Brownian motion started
at $x$, and $\BMM_\iv$ is its local time at $0$ in the normalization of Skorokhod's reflection problem. That is,
$\BMM$ is continuous and nondecreasing, it increases only when $\BMr=0$, and $\BMr=\BMM-\BM$ stays nonnegative. See
\citet{skorokhod1961stochastic} and \citet{levyProcessusStochastiquesMouvement1965} for the classical results, and \citet[Ch.~VI, Lemma~(2.1) and Thm.~(2.3)]{RevuzYor99} for a
modern account. In the convention of Revuz and Yor, in which the local time $L^0(\BMr)$ is defined by Tanaka's
formula, $\BMM_t=\frac12L^0_t(\BMr)$ \citep[Ch.~VI, Exercise~(1.16)]{RevuzYor99}. Figure~\ref{fig:reflection}
shows this construction of $\BMr$.

\begin{figure}[htb]
	\centering
	\includegraphics[width=0.8\textwidth]{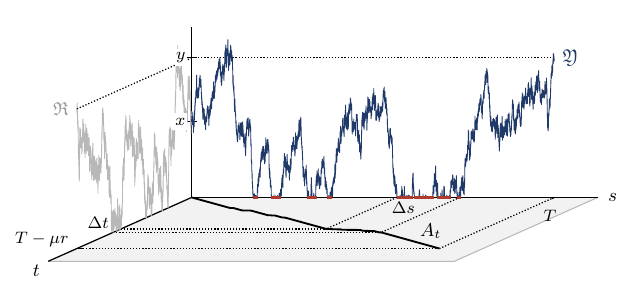}
	\caption[The time change]{The path $\BMr$ lies in the $(t,\text{position})$ plane, the path $\Y$ in the
	$(s,\text{position})$ plane, and the clock $A_t=t+\stick\BMM_t$ in the $(t,s)$ plane. Both paths start at $x$ and
	end at $y=\Y_T=\BMr_{T-\stick r}$, where $r=\locT_T$. The zero set of $\Y$ over time $T$ has
	Lebesgue measure $\stick r$. The inverse clock $\sigma$ takes an interval of length $\Delta s$ to one of length 
	$\Delta t=\Delta s-\stick\Delta\locT$, where $\Delta\locT$ is the local time that $\Y$ gains over the interval. 
	The path is that of Figure~\ref{fig:reflection}, with $\stick=0.6$ and $T=3.3$.}
	\label{fig:R_A_Y}
\end{figure}

Given $\stick>0$, the sticky clock is $A_t\ass t+\stick\BMM_t$. It is continuous and strictly increasing, and we write
$\sigma_\iv$ for its inverse. Sticky reflected Brownian motion and its local time are
\[
	\Y_s~\ass~\BMr_{\sigma_s},\qquad \locT_s~\ass~\BMM_{\sigma_s},\qquad s\geq0.
\]
Figure~\ref{fig:R_A_Y} shows the path of Figure~\ref{fig:reflection} on both clocks at once.
We call $t$ the operational time and $s=A_t$ the physical time. Every unit of local time adds $\stick$ units of physical
time at $0$. The inverse clock $\sigma_s=s-\stick\locT_s$ is continuous and strictly increasing. The zero set of
$\BMr$ has Lebesgue measure zero, and therefore
\[
	\int_0^s \id_{\{0\}}(\Y_u)\md u~=~\stick\,\locT_s .
\]
The local time increases only when $\Y=0$, so this identity gives $\md\locT_s=\stick^{-1}\id_{\{0\}}(\Y_s)\md s$. The
inverse clock is therefore absolutely continuous, with derivative $\id_{(0,\infty)}(\Y_s)$ at almost every $s$.
The local time of $\Y$ reaches $a$ at the physical time $A_{\hT_a}=\hT_a+\stick a$. The excursion intervals of $\Y$ are
those of $\BMr$, each shifted by $\stick$ times the local time at the start of the interval. We
take this time change as the definition of sticky reflected Brownian motion. 

In Feller's classification, $\Y$ is the Feller diffusion on
$[0,\infty)$ with generator $\frac12\pd^2_{xx}$ on functions satisfying $f'(0)=\frac{\stick}{2}f''(0)$
\citep{fellerParabolicDifferentialEquations1952,itoBrownianMotionsHalfline1963}.
The process also solves the stochastic differential equation
$\md\Y_s=\id_{(0,\infty)}(\Y_s)\md\beta_s+\stick^{-1}\id_{\{0\}}(\Y_s)\md s$ driven by a Brownian motion $\beta$, which has a weak solution
that is unique in law but no strong solution \citep{chitashviliNonexistenceStrongSolution1989}.
\citet{warrenBranchingProcessesRayKnight1997} found the conditional law of the position given $\beta$. The time change above
builds the path from $\BM$, which is a different Brownian motion from $\beta$.
\FloatBarrier

\subsection{Joint transition law of position and local time}\label{sec:transition}

The samplers rest on the joint law of the position and the local time at a fixed time. We derive this law for
reflected Brownian motion from the reflection principle and transfer it through the time change. For $t>0$ we write
\[
	\heat(a,t)~\ass~\frac{\me^{-a^2/2t}}{\sqrt{2\mpi t}},\qquad h(a,t)~\ass~\frac at\,\heat(a,t),\qquad
	p_{\mathfrak A}(y,t\mid x)~\ass~\heat(y-x,t)-\heat(y+x,t)
\]
for the heat kernel, the density at time $t$ of the first passage of Brownian motion over a distance $a>0$
(Section~\ref{sec:passage}), and the transition density of Brownian motion killed at $0$.

\begin{proposition}\label{prop:rbm-joint}
For $x\geq0$, $t>0$, $y\geq0$ and $r>0$,
\begin{align*}
	\Probx{x}\{\BMr_t\in\md y,\ \BMM_t=0\}~&=~p_{\mathfrak A}(y,t\mid x)\md y,\\
	\Probx{x}\{\BMr_t\in\md y,\ \BMM_t\in\md r\}~&=~2h(x+y+r,t)\md y\md r .
\end{align*}
\end{proposition}

\begin{proof}
On $\{M_t\geq0\}$, with $M_t\ass\max_{u\leq t}\BM_u$, the pair $(\BMr_t,\BMM_t)$ is $(M_t-\BM_t,M_t)$, and the joint
density $2h(2m-b+x,t)$ of $(M_t,\BM_t)$ from the reflection principle (\citealp{levyProcessusStochastiquesMouvement1965};\ \citealp[Prop.~2.8.1]{karatzas1991brownian}) gives the second line. On
$\{M_t<0\}$ the path has not reached $0$, and $\BMr_t=-\BM_t$ is Brownian motion killed at $0$.
\end{proof}

We now transfer the joint law of $(\BMr,\BMM)$ through the sticky clock. We write
$\pYl(y,r,s\mid x)$ for the joint transition density of position and local time,
starting from $(x,0)$, defined by
\[
	\Probx{x}\{\Y_s\in\md y,\ \locT_s\in\md r\}
	~=~\pYl(y,r,s\mid x)
	\big(\delta_0(\md y)+\md y\big)
	\big(\delta_0(\md r)+\md r\big).
\]
The following joint law is given by
\citet[Thm.~2]{casterasTrivariateDistributionSticky2026} in a different
normalization. We derive it here directly from the time change.

\begin{proposition}\label{prop:kernel}
For $x\geq0$, $s>0$ and $y,r\geq0$,
\[
	\pYl(y,r,s\mid x)
	~=~
	\begin{cases}
		p_{\mathfrak A}(y,s\mid x),
			& y>0,\ r=0,\\[2pt]
		2h(x+y+r,\ s-\stick r),
			& y>0,\ 0<r<s/\stick,\\[2pt]
		\stick\,h(x+r,\ s-\stick r),
			& y=0,\ 0<r<s/\stick,\\[2pt]
		0,
			& r\ge s/\stick, \text{ or } y=0 \text{ and } r=0.
	\end{cases}
\]
\end{proposition}

\begin{proof}
Let $F\geq0$ be continuous with compact support in
$(0,\infty)\times[0,\infty)^2$. The clock $A$ is continuous and strictly
increasing, with
\[
	\md A_t=\md t+\stick\md\BMM_t.
\]
The substitution $s=A_t$ gives
\[
	\E_{x}\int_0^\infty F(s,\Y_s,\locT_s)\md s
	~=~
	\E_{x}\int_0^\infty F(A_t,\BMr_t,\BMM_t)\md t
	+\stick\,\E_{x}\int_0^\infty F(A_t,0,\BMM_t)\md\BMM_t,
\]
where the second term uses that $\BMM$ increases only when $\BMr=0$.

The set $\{t:\BMr_t=0\}$ has Lebesgue measure zero, so the first term
puts no mass at $\Y=0$. Proposition~\ref{prop:rbm-joint} and the
substitution $s=t+\stick r$, with Jacobian $1$, give the first two
cases of the density.

For the second term, we change variables from operational time to the
local time level $r$. The process $\BMM$ reaches $r$ at $\hT_r$, whose
density under $\Probx{x}$ is $h(x+r,\iv)$. Hence
\begin{align*}
	\stick\,\E_{x}\int_0^\infty F(A_t,0,\BMM_t)\md\BMM_t
	~&=~
	\stick\int_0^\infty
	\E_{x}\big[F(\hT_r+\stick r,0,r)\big]\md r\\
	&=~
	\int_0^\infty\int_0^{s/\stick}
	F(s,0,r)\,\stick\,h(x+r,\ s-\stick r)\md r\md s.
\end{align*}
This gives the third case. Thus the local time term
$\stick\md\BMM_t$ of the clock supplies the mass at $\Y_s=0$.

Finally, take $F(s,y,r)=\phi(s)\,G(y,r)$, where $\phi$ and $G$ are
nonnegative continuous functions with compact support. The function
\[
	s\longmapsto \E_x\big[G(\Y_s,\locT_s)\big]
\]
is continuous by continuity of the paths and bounded convergence.
The corresponding expression obtained from the displayed densities
is also continuous for $s>0$ by dominated convergence, since the
densities are locally bounded in $s$ on the support of $G$.
Equality of their integrals against every $\phi$ therefore gives
equality at every $s>0$, proving the stated joint law.
\end{proof}

The support restriction $r<s/\stick$ follows from
$\sigma_s=s-\stick\locT_s>0$ for $s>0$. Integrating the joint density
over the local time gives
\[
	\Probx{x}\{\Y_s\in\md y\}
	~=~
	\left(
		\pYl(y,0,s\mid x)
		+\int_0^{s/\stick}\pYl(y,r,s\mid x)\md r
	\right)
	\big(\delta_0(\md y)+\md y\big).
\]
This agrees with the transition kernel in
\citet[eq.~(3.1)]{hajriApplicationStochasticFlows2017}, with their
parameter $\theta=1/\stick$. They obtain that kernel by inverse
Laplace transformation of the resolvent given in
\citet[Prop.~13]{warrenBranchingProcessesRayKnight1997}.

\subsection{Strong Markov property of the joint process}\label{sec:markov}

The process $\Y$ is a Feller diffusion and is therefore strong Markov.
The identity
\[
	\locT_{s+v}-\locT_s
	~=~\stick^{-1}\int_0^v\id_{\{0\}}(\Y_{s+w})\md w,
	\qquad s,v\geq0,
\]
shows that the increment of the local time is a function of the path of $\Y$ after $s$. The strong Markov property
therefore extends to the pair $(\Y,\locT)$, with the local time already gained entering only as a shift.

For the calculations below, we use the following consequence of the
strong Markov property of $\BM$, expressed in terms of the filtration
$(\filt_t)$ of Section~\ref{sec:time-change}. Let $S$ be an almost surely
finite nonnegative random time such that $\sigma_S$ is a stopping time
of $\BM$. For every Borel set $\Gamma$ of continuous paths with values
in $[0,\infty)^2$,
\begin{equation}\label{eq:markov}
	\Probx{x}\big\{
		(\Y_{S+\iv},\,\locT_{S+\iv}-\locT_S)\in\Gamma
		\bigm|\filt_{\sigma_S}
	\big\}
	~=~
	\Probx{\Y_S}\big\{
		(\Y_\iv,\,\locT_\iv)\in\Gamma
	\big\},
	\qquad \Probx{x}\text{-a.s.}
\end{equation}
Thus the conditional law of the future position and local time
increment depends on the past only through $\Y_S$.
Lemma~\ref{lem:markov} derives \eqref{eq:markov} from the strong
Markov property of $\BM$ at $\sigma_S$.

For $u\geq0$, define
\[
	g_u~\ass~\sup\{v\leq u:\Y_v=0\},
	\qquad
	d_u~\ass~\inf\{v\geq u:\Y_v=0\},
\]
with $g_u\ass0$ if there is no zero in $[0,u]$ and
$\inf\varnothing\ass\infty$. These are the last zero at or before
$u$ and the first zero at or after $u$. In particular,
$d_0=\hT_0$, since no local time is gained before the first zero
and the physical and operational times agree up to it.

Lemma~\ref{lem:markov} verifies the stopping-time condition at fixed
times, at $d_u$ for fixed $u$, and at $\hT_a+\stick a$ for fixed
$a>0$. It also verifies the condition at a finite
$\filt_{\sigma_S}$-measurable time $S'\geq S$ whenever it holds
at $S$. Consequently, \eqref{eq:markov} applies successively at
times $s_0<s_1<\cdots$ when $s_0$ is deterministic and each
$s_{k+1}$ is an almost surely finite measurable function of
$(s_j,\Y_{s_j},\locT_{s_j})_{j\leq k}$.
The forward sampler uses \eqref{eq:markov} at these successive times.

Inserting a value between previously sampled times requires its
conditional law given the retained observations on both sides,
including any recorded zeros. The identity \eqref{eq:markov} alone
does not supply this law. Theorem~\ref{thm:refine} treats these queries
through the conditional law of the path given the skeleton
(Appendix~\ref{app:path-record}).

For conditioning on specified endpoint values, including
$\{\Y_T=y,\ \locT_T=r\}$, we use the versions defined by the
density ratios in Appendix~\ref{app:conditional}.

\subsection{Conditioned passage times and Bessel(3) bridges}\label{sec:passage}

Brownian passage times will be used to sample zeros and local time, while Bessel(3) bridges will supply values on
intervals without zeros. Appendix~\ref{app:brownian} gives the samplers for conditioned passage times and Bessel(3)
bridges and proves their laws. We first work under $\Probx{0}$, so that $\BM$ starts at $0$ and $\hT_a$ is its first
passage time to level $a$.

For $a>0$, this passage time has the density
\begin{equation}\label{eq:h}
	\Probx{0}\big\{\hT_a\in\md t\big\}~=~h(a,t)\md t,\qquad t>0,
\end{equation}
where $h$ is the function of Section~\ref{sec:transition},
\[
	h(a,t)~=~\frac{a}{t}\,\heat(a,t)~=~\frac{a}{\sqrt{2\mpi t^3}}\,\me^{-a^2/2t}.
\]
Under $\Probx{x}$, the corresponding density is $h(a+x,t)$. We set $h(a,t)\ass0$ for $t\leq0$.
By the reflection principle, for $t>0$,
\[
	\Probx{0}\{\hT_a\leq t\}
	~=~\Probx{0}\{|\BM_t|\geq a\}
	~=~\Prob\{a^2/\gauss^2\leq t\},
\]
where $\gauss$ is standard normal. Hence
\begin{equation}\label{eq:h-sampler}
	\hT_a~\law~a^2/\gauss^2 .
\end{equation}
Under $\Probx{0}$, the strong Markov property at $\hT_a$ implies that
$\hT_{a+b}-\hT_a$ is independent of $\hT_a$ and has the law of $\hT_b$, for $a,b>0$. Therefore
\begin{equation}\label{eq:h-conv}
	\int_0^t h(a,u)\,h(b,t-u)\md u~=~h(a+b,t),\qquad t>0 .
\end{equation}

The passage time density also satisfies an identity with the transition density $p_{\mathfrak A}$ of Brownian
motion killed at $0$. Section~\ref{sec:bridge} uses this identity to derive the distribution of the last zero
before a fixed time,
\begin{equation}\label{eq:entrance}
	\int_0^\infty h(w,v)\,p_{\mathfrak A}(y,t\mid w)\md w~=~h(y,v+t),\qquad v,t,y>0 .
\end{equation}
Thus the measures $h(y,v)\md y$, $v>0$, form an entrance law for Brownian motion killed at $0$.
To verify \eqref{eq:entrance}, extend $h(\cdot,v)$ to an odd function on $\R$. Since
$p_{\mathfrak A}(y,t\mid w)=\heat(y-w,t)-\heat(y+w,t)$ and
$h(w,v)=-\pd_w\heat(w,v)$, integration by parts and the Gaussian convolution identity give
\begin{align*}
	\int_0^\infty h(w,v)\,p_{\mathfrak A}(y,t\mid w)\md w
	&= \int_{\R}h(w,v)\,\heat(y-w,t)\md w\\
	&= -\pd_y\int_{\R}\heat(w,v)\,\heat(y-w,t)\md w\\
	&= -\pd_y\heat(y,v+t)
	~=~h(y,v+t).
\end{align*}

The bridge construction requires a passage time conditioned on a later passage time. For $a,b,t>0$, the joint
density of $\hT_a$ and $\hT_{a+b}$ at $(u,t)$ is $h(a,u)h(b,t-u)$ on $0<u<t$. Taking the ratio with the density of
$\hT_{a+b}$ gives
\begin{equation}\label{eq:ptb}
	\Probx{0}\big\{\hT_a\in\md u\mid\hT_{a+b}=t\big\}
	~=~\frac{h(a,u)\,h(b,t-u)}{h(a+b,t)}\md u,\qquad 0<u<t,
\end{equation}
which integrates to one by \eqref{eq:h-conv}. Proposition~\ref{prop:ptb} represents this law as a mixture of two
transformed inverse Gaussian distributions, the passage time form of the sampler for the time of the minimum in
\citet[Prop.~2]{beskosRetrospectiveExactSimulation2006}. Selecting the component first and generating only the
selected inverse Gaussian variable uses one normal and two uniforms, without rejection.

Let $F_{m,\lambda}$ denote the distribution function of the inverse Gaussian law $\IG(m,\lambda)$ with mean $m$
and shape $\lambda$. The mixture representation also gives
\begin{equation}\label{eq:ptb-cdf}
	\begin{aligned}
		\Probx{0}\big\{\hT_a\leq u\mid\hT_{a+b}=t\big\}
		~&=~\frac{b}{a+b}\,
		F_{a/b,\,a^2/t}\Big(\frac{u}{t-u}\Big)\\
		&\quad+\frac{a}{a+b}
		\Big(1-F_{b/a,\,b^2/t}\Big(\frac{t}{u}-1\Big)\Big),
		\qquad 0<u<t .
	\end{aligned}
\end{equation}

The Bessel(3) process is the Doob $\mathfrak h$-transform of Brownian motion killed at $0$, with
$\mathfrak h(x)=x$
\citep{doobConditionalBrownianMotion1957a,williamsPathDecompositionContinuity1974,pitmanOnedimensionalBrownianMotion1975}.
For $z,z'>0$ and $t>0$, its bridge from $z$ to $z'$ over $[0,t]$ coincides with Brownian motion started at $z$,
conditioned to stay positive up to time $t$ and to end at $z'$. Appendix~\ref{app:bridges-paths} defines these
bridges, including those with an endpoint at $0$, as laws of paths through a three-dimensional Brownian motion.

Write $\bessb$ for the Bessel(3) bridge from $z$ to $z'$ of duration $t$.
The samplers require its value $\bessb_u$ at $0<u<t$ in three cases,
\[
	z=0<z',\qquad z'=0<z,\qquad z,z'>0 .
\]
The first describes the path after its last zero, the second the path before its first zero, and the third the path
between two positive values when the local time remains constant.

A Bessel(3) bridge can be constructed as the length of a three-dimensional Brownian motion conditioned on the length
of its endpoint. Conditional on the endpoint direction, its value at time $u$ is the length of a Gaussian vector.
When both endpoints are positive, the law of this length depends on the direction through the cosine of the angle
between the initial and terminal vectors, which must first be sampled. When either endpoint is $0$, the law of the
length does not depend on this direction. Formula~\eqref{eq:bessel-sampler} of Appendix~\ref{app:brownian} samples
the bridge value using one standard normal, one exponential and, when both endpoints are positive, one uniform.
Lemma~\ref{lem:bessel} proves its law.

\section{Forward sampling of position and local time}\label{sec:forward}

We now construct the forward step for the pair $(\Y,\locT)$. Starting from $x\geq0$, we draw
$(\Y_s,\locT_s)$ under $\Probx{x}$ from the transition law of Proposition~\ref{prop:kernel}. When the step starts from a
state $(x,a)$, the local time returned is the increment, which is added to $a$. We sample this increment first. Its tail
probabilities are Brownian passage probabilities, and a single half-normal variable determines whether the path reaches
$0$ and how much local time it gains. Once one has the local time, the position needs at most one uniform and one
exponential.

\subsection{Local time distribution and inversion}\label{sec:fwd-lt}

For a fixed $u>0$, the local time first reaches $u$ at $\hT_u+\stick u$
(Section~\ref{sec:time-change}), so $\locT_s\geq u$ exactly when $\hT_u+\stick u\leq s$. The two events
$\{\locT_s>u\}$ and $\{\hT_u+\stick u<s\}$ differ only on $\{\locT_s=u\}$, which has probability zero. Indeed, by
the strong Markov property at $\hT_u$ and the reflection principle, $\BM$ exceeds $u$ immediately after $\hT_u$,
almost surely, so the local time stays at the level $u$ for no positive time. Thus, for $0<u<s/\stick$,
\begin{align*}
	\Probx{x}\{\locT_s>u\}~&=~\Probx{x}\{\hT_u+\stick u<s\}\\
	&=~\Prob\big\{(x+u)^2/\gauss^2<s-\stick u\big\}\\
	&=~\Prob\big\{|\gauss|>(x+u)/\sqrt{s-\stick u}\big\}.
\end{align*}
The second line uses \eqref{eq:h-sampler}, since $\hT_u$ is the first passage of $\BM$ from $-x$ to $u$.
Letting $u\downarrow0$ gives
\[
	\Probx{x}\{\locT_s=0\}~=~\Prob\big\{|\gauss|\leq x/\sqrt s\big\}.
\]

The function
\[
	u~\mapsto~\frac{x+u}{\sqrt{s-\stick u}}
\]
is continuous and strictly increasing from $x/\sqrt s$ to $\infty$ on $[0,s/\stick)$. Hence we draw
$H\ass|\gauss|$ and return local time $0$ if $H\leq x/\sqrt s$. Otherwise we take the unique
$r\in(0,s/\stick)$ satisfying
\[
	H\sqrt{s-\stick r}~=~x+r.
\]
Squaring gives
\[
	r^2+(2x+\stick H^2)\,r+x^2-sH^2~=~0.
\]
Rationalizing the positive root gives
\begin{equation}\label{eq:fwd-root}
	r~=~\frac{2(sH^2-x^2)}
	{2x+\stick H^2+H\sqrt{\stick^2H^2+4(s+\stick x)}}.
\end{equation}

\subsection{The position conditioned on the local time}\label{sec:fwd-pos}

For $0<r<s/\stick$, put $c\ass x+r$ and $t\ass s-\stick r$. By Proposition~\ref{prop:kernel}, the conditional
law of $\Y_s$ given $\locT_s=r$ is
\[
	\frac{\stick h(c,t)\,\delta_0(\md y)
		+2h(c+y,t)\,\id_{(0,\infty)}(y)\md y}
	{\stick h(c,t)+2\heat(c,t)},
\]
since $\int_0^\infty2h(c+y,t)\md y=2\heat(c,t)$. Using
$h(c,t)=\frac ct\heat(c,t)$, we obtain
\begin{align*}
	\Probx{x}\{\Y_s=0\mid\locT_s=r\}
	~&=~\frac{\stick c}{\stick c+2t},\\
	\Probx{x}\{\Y_s>y\mid\locT_s=r,\ \Y_s>0\}
	~&=~\frac{\int_y^\infty2h(c+z,t)\md z}{2\heat(c,t)}
	~=~\exp\Big\{-\frac{(c+y)^2-c^2}{2t}\Big\},
	\qquad y\geq0.
\end{align*}
Conditional on $\Y_s>0$, the position therefore has the law of
$\sqrt{c^2+2t\erv}-c$, where $\erv\law\Exp(1)$.

For $x>0$, conditional on $\locT_s=0$, the path has not reached $0$ by time $s$, almost surely. Physical and
operational time then agree, and $\Y_s=-\BM_s$. Write $\W\ass\BM+x$, so that $\W_0=0$, and put
$M_s\ass\max_{u\leq s}\W_u$. The reflection principle gives $M_s\law\sqrt s\,H$. Consequently, conditional on
$H\leq x/\sqrt s$, the variable $m\ass\sqrt s\,H$ has the law of $M_s$ conditional on $M_s\leq x$.
We reuse $H$ to draw this maximum and then draw the endpoint given it.

For $m>0$, the conditional endpoint law is
\begin{equation}\label{eq:seshadri}
	\W_s\mid\{M_s=m\}~\law~2m-\sqrt{m^2+2s\erv},
	\qquad \erv\law\Exp(1).
\end{equation}
Indeed, the reflection principle \citep{levyProcessusStochastiquesMouvement1965} gives
\[
	\Probx{x}\{2M_s-\W_s>v\mid M_s=m\}
	~=~\frac{\heat(v,s)}{\heat(m,s)}
	~=~\exp\Big\{-\frac{v^2-m^2}{2s}\Big\},
	\qquad v\geq m,
\]
which yields \eqref{eq:seshadri}.\footnote{\citet[Sec.~6]{seshadriExponentialModelsBrownian1988a} shows that
$M_s(M_s-\W_s)$ is independent of $\W_s$, which gives the maximum given the endpoint. The law
\eqref{eq:seshadri} is the reverse conditioning, the endpoint given the maximum.}
The returned position is $x-\W_s$. When $x=0$, the event $\{\locT_s=0\}$ has probability zero.

\subsection{The forward algorithm and repeated sampling}\label{sec:fwd-alg}

Algorithm~\ref{alg:forward} combines the two previous subsections. A step uses one normal, at most one uniform and at
most one exponential, and no rejection. Here $\gauss\law N(0,1)$, $\uni\law\mathrm{Unif}(0,1)$ and
$\erv\law\Exp(1)$, and each draw is independent of all preceding draws.

\begin{algorithm}[H]
\caption{The forward step (Theorem~\ref{thm:forward})}\label{alg:forward}
\begin{algorithmic}[1]
\Require a value $x\geq0$ and a time $s>0$
\State draw $\gauss$ and put $H\ass|\gauss|$
\If{$H\leq x/\sqrt s$}\Comment{the path does not reach $0$, almost surely}
	\State draw $\erv$, put $m\ass\sqrt s\,H$, and \Return
	$\big(x-2m+\sqrt{m^2+2s\erv},\ 0\big)$\Comment{\eqref{eq:seshadri}}
\EndIf
\State set $r$ by \eqref{eq:fwd-root}, and put $c\ass x+r$ and $t\ass c^2/H^2=s-\stick r$
\State draw $\uni$; \textbf{if} $\uni\leq\stick c/(\stick c+2t)$ \Return $(0,r)$
\State draw $\erv$ and \Return $\big(\sqrt{c^2+2t\erv}-c,\ r\big)$
\end{algorithmic}
\end{algorithm}

The repeated forward step starts from $s_0=0$ and $(Y_0,L_0)\ass(x,0)$. At step $i$, it chooses a finite time
$s_i>s_{i-1}$, draws $(Y_i,R_i)$ by Algorithm~\ref{alg:forward} from $Y_{i-1}$ over the time $s_i-s_{i-1}$,
and sets $L_i\ass L_{i-1}+R_i$. Each step uses fresh variables. The times may be fixed in advance or chosen by
measurable functions $f_i$, so that
\[
	s_i~=~f_i\big((Y_j,L_j)_{j<i}\big).
\]
We require the $f_i$ to give finite, strictly increasing times. On the path, define
$S_0\ass0$ and
\[
	S_i~\ass~f_i\big((\Y_{S_j},\locT_{S_j})_{j<i}\big).
\]
Each $S_{i+1}$ is measurable with respect to $\filt_{\sigma_{S_i}}$ and exceeds $S_i$, so
Lemma~\ref{lem:markov} applies successively at the times $S_i$.

\begin{theorem}\label{thm:forward}
For $x\geq0$ and $s>0$, the pair $(Y,L)$ returned by Algorithm~\ref{alg:forward} has the law of $(\Y_s,\locT_s)$
under $\Probx{x}$, that is,
\[
	\Prob\{Y\in\md y,\ L\in\md r\}
	~=~\pYl(y,r,s\mid x)
	\big(\delta_0(\md y)+\md y\big)
	\big(\delta_0(\md r)+\md r\big).
\]
For the repeated forward step and the times defined above,
\[
	\big(s_i,\,Y_i,\,L_i\big)_{i\leq n}
	~\law~
	\big(S_i,\,\Y_{S_i},\,\locT_{S_i}\big)_{i\leq n},
\]
where the process on the right has law $\Probx{x}$.
\end{theorem}

\begin{proof}
Let $L$ be the local time returned by the algorithm. The monotonicity of the threshold in
Section~\ref{sec:fwd-lt} gives, for $0\leq u<s/\stick$,
\[
	\Prob\{L>u\}
	~=~\Prob\big\{H>(x+u)/\sqrt{s-\stick u}\big\}
	~=~\Probx{x}\{\locT_s>u\}.
\]
Both variables take values in $[0,s/\stick)$, so $L$ has the law of $\locT_s$.

For $0<r<s/\stick$, the algorithm draws the atom at $0$ and the positive position from the conditional law
in Section~\ref{sec:fwd-pos}. The uniform and exponential used for these draws are independent of $H$, and hence of $L$.

On the branch $L=0$, the algorithm draws an endpoint given the maximum
$m=\sqrt s\,H\leq x$ by \eqref{eq:seshadri} and returns $x-\W_s$. The joint density of $(M_s,\W_s)$ is
$2h(2m-w,s)$ for $m>\max(0,w)$. Thus, for $y>0$,
\begin{align*}
	\Prob\{Y\in\md y,\ L=0\}
	~&=~
	\left[\int_{\max(0,\,x-y)}^x2h(2m-x+y,s)\md m\right]\md y\\
	&=~\big[\heat(x-y,s)-\heat(x+y,s)\big]\md y\\
	&=~p_{\mathfrak A}(y,s\mid x)\md y.
\end{align*}
This is the part of the transition law in Proposition~\ref{prop:kernel} with zero local time, and proves the claim
for one step.

For repeated steps, suppose the sampled triples have the stated joint law through step $i$. On the process side,
$S_{i+1}$ is determined by the observations through $S_i$. By Lemma~\ref{lem:markov}, conditional on those
observations,
\[
	\big(\Y_{S_{i+1}},\,\locT_{S_{i+1}}-\locT_{S_i}\big)
\]
has the transition law from $\Y_{S_i}$ over the known duration $S_{i+1}-S_i$. On the algorithm side, the next time
is chosen by the same function, and the next step draws from that transition law with fresh variables. Adding the local
time increment gives the same conditional law for the next triple. Induction on $i$ proves the result.
\end{proof}

\section{Bridges conditioned on position and local time}\label{sec:bridge}

In this section we sample $(\Y_s,\locT_s)$ at a time $s\in(0,T)$, conditional on $(\Y_T,\locT_T)=(y,r)$. Fix
$T>0$, $x,y\geq0$ and $r>0$ with $\stick r<T$. We first draw the first and last zeros of the bridge from
conditioned Brownian passage laws. Outside these zeros, the value comes from a Bessel(3) bridge. Between them, we draw
the local time by inversion and then the position by rejection, in the order of Section~\ref{sec:forward}. All of the
local time $r$, and all of the time $\stick r$ spent at $0$, lie between the two zeros, since the local time is
constant while the path is away from $0$. We keep the zeros together with the sampled value, because the refinement of
Section~\ref{sec:refine} needs them. The case $r=0$, in which the path has no zero, closes
Section~\ref{sec:bridge-alg}.

The bridge is the conditional law \eqref{eq:bridge-def} given $(\Y_T,\locT_T)=(y,r)$. It is defined at every endpoint with
$\pYl(y,r,T\mid x)>0$. Conditional laws given further random variables, such as the zeros below, are ratios of densities
in the same way. For $V=(\Y_s,\locT_s)$, the Markov property at $s$ (Lemma~\ref{lem:markov}) gives the law that the
section samples,
\begin{equation}\label{eq:bridge-one}
	\begin{aligned}
		\Probx{x}\{\Y_s\in\md z,\ \locT_s\in\md a\mid\Y_T=y,\ \locT_T=r\}
		~&=~\frac{\pYl(z,a,s\mid x)\,\pYl(y,r-a,T-s\mid z)}{\pYl(y,r,T\mid x)}\\
		&\qquad\cdot\big(\delta_0(\md z)+\md z\big)\big(\delta_0(\md a)+\md a+\delta_r(\md a)\big)
	\end{aligned}
\end{equation}
on $[0,\infty)\times[0,r]$. Since $r>0$, the two endpoint atoms are
distinct. At $a=0$ the path gains no local time on $[0,s]$; at $a=r$
it gains no local time on $[s,T]$. When
$y>0$ and $\stick r<s$, the bridge can gain all of its local time before $s$ and stay positive from $g_T<s$ to $T$, so
$\{\locT_s=r\}$ has positive probability. The local time $a$ enters \eqref{eq:bridge-one} through the positions and
through the operational times $s-\stick a$ and $T-s-\stick(r-a)$ inside the two factors. We know no closed form for
$\Probx{x}\{\Y_s=0\mid\Y_T=y,\ \locT_T=r\}$ or for the distribution function of either coordinate, so we do not
sample \eqref{eq:bridge-one} directly. The zeros and the three bridges between them have explicit laws, and this is
what the decomposition gives.

\subsection{The first and last zeros}\label{sec:bridge-zeros}

Recall from Section~\ref{sec:markov} that $d_u$ is the first zero of
$\Y$ at or after $u$, and $g_u$ is its last zero at or before $u$,
with $g_u=0$ if the path has yet to hit zero. In particular,
$d_0=\hT_0$ is the first zero of the path. Suppose $\locT_u=a>0$ and $\Y_u>0$.
The last zero of $\BMr$ before operational time $\sigma_u$ occurs at $\hT_a$, when $\BM$ first reaches $a$. Hence
\begin{equation}\label{eq:last-zero}
	g_u~=~A_{\hT_a}~=~\hT_a+\stick a\qquad\text{on }\{\locT_u=a,\ \Y_u>0\},
\end{equation}
and $g_u=u$ on $\{\Y_u=0\}$. We now find the joint law of the first zero, the last zero and the state at a fixed
time. Divided by the density of the state, it gives the two draws of the zeros in Corollary~\ref{cor:zeros}. The
computation rests on the entrance identity \eqref{eq:entrance}.

\begin{proposition}\label{prop:zeros}
Let $x>0$, $s>0$ and $0<r<s/\stick$. For $y>0$, $f>0$ and $f+\stick r<g<s$,
\begin{align*}
	\Probx{x}\{\hT_0\in\md f,\ g_s\in\md g,\ \Y_s\in\md y,\ \locT_s\in\md r\}
	~&=~h(x,f)\,h(r,\ g-f-\stick r)\,2h(y,\ s-g)\md f\md g\md y\md r,\\
	\Probx{x}\{\hT_0\in\md f,\ \Y_s=0,\ \locT_s\in\md r\}~&=~h(x,f)\,\stick\,h(r,\ s-f-\stick r)\md f\md r,
\end{align*}
the second for $0<f<s-\stick r$. Outside these ranges the densities vanish, by the convention $h(a,t)=0$ for $t\leq0$.
\end{proposition}

On $\{\Y_s=0\}$ we have $g_s=s$. For $x=0$ the same holds with $\hT_0=0$, and with the point mass at $f=0$ in place of
$h(x,f)\md f$.

\begin{proof}
The local time is $0$ before $\hT_0$, so the physical and operational times agree there. Hence $\hT_0$ is the passage
time of $\BM$ from $-x$ to $0$, with density $h(x,\iv)$ by \eqref{eq:h}. By Lemma~\ref{lem:markov} at $S=\hT_0$, where
$(\Y,\locT)=(0,0)$,
\[
	\Probx{x}\big\{\hT_0\in\md f,\ (g_s,\Y_s,\locT_s)\in\iv\big\}~=~h(x,f)\md f\;
	\Probx{0}\big\{(f+g_{s-f},\ \Y_{s-f},\ \locT_{s-f})\in\iv\big\},
\]
and it is enough to take $x=0$. On $\{\Y_s=0\}$ the formula is the third line of Proposition~\ref{prop:kernel}, and
$g_s=s$ by definition. Let $y>0$ and $\stick r<g<s$. The last zero is at most $g$ exactly when $\Y$ has no zero in
$(g,s]$. From a state $(w,r)$ with $w>0$, the path up to its first zero is a Brownian motion from $w$ with constant local
time. From $w=0$, the path has a zero at once. The Markov property at the fixed time $g$, Proposition~\ref{prop:kernel},
\eqref{eq:h-conv} and \eqref{eq:entrance} give
\begin{align*}
	\Probx{0}\{g_s\leq g,\ \Y_s\in\md y,\ \locT_s\in\md r\}
	~&=~\int_{w>0}\Probx{0}\{\Y_g\in\md w,\ \locT_g\in\md r\}\,p_{\mathfrak A}(y,s-g\mid w)\md y\\
	&=~\int_0^\infty2h(w+r,\ g-\stick r)\,p_{\mathfrak A}(y,s-g\mid w)\md w\md y\md r\\
	&=~2\int_0^{g-\stick r}h(r,v)\int_0^\infty h(w,\ g-\stick r-v)\,p_{\mathfrak A}(y,s-g\mid w)\md w\md v\\
	&\qquad\qquad\cdot\md y\md r\\
	&=~2\int_0^{g-\stick r}h(r,v)\,h(y,\ s-\stick r-v)\md v\md y\md r .
\end{align*}
Differentiating in $g$ gives the density $2h(r,g-\stick r)\,h(y,s-g)$, which is the first line for $x=0$.
\end{proof}

After the first zero, and in operational time, Proposition~\ref{prop:zeros} is the joint law of the maximum $r$ of
a Brownian path on $[0,s-\stick r]$, its location $g_s-\stick r$ by \eqref{eq:last-zero}, and the endpoint
\citep{levyProcessusStochastiquesMouvement1965,sheppJointDensityMaximum1979}. \citet{imhofDensityFactorizationsBrownian1984}
wrote this law as the product of the two passage-time densities $h(r,\iv)$ and $2h(y,\iv)$ of the first line, and it is
the decomposition of a Brownian path at its maximum
\citep{williamsPathDecompositionContinuity1974,pitmanDecompositionMaximumExcursions1996}; see
\citet[Sec.~2.8]{karatzas1991brownian}. The factor $h(x,f)$ is the first zero, and the shift by $\stick r$ is the time
change. Integrating over $f$ and $g$ by \eqref{eq:h-conv} gives back Proposition~\ref{prop:kernel}, and dividing by it
gives the zeros of the bridge. When $x=0$ the first zero is $\hT_0=0$ and needs no draw. When $y=0$ the last zero is
$T$, and $\hT_r-\hT_0=T-\stick r-\hT_0$ needs no draw either.

\begin{corollary}\label{cor:zeros}
Let $y>0$. Then
\begin{align*}
	\Probx{x}\{\hT_0\in\md f\mid\Y_T=y,\ \locT_T=r\}
	~&=~\Probx{0}\{\hT_x\in\md f\mid\hT_{x+r+y}=T-\stick r\},\\
	\Probx{x}\{\hT_r-\hT_0\in\md t\mid\hT_0=f,\ \Y_T=y,\ \locT_T=r\}
	~&=~\Probx{0}\{\hT_r\in\md t\mid\hT_{r+y}=T-\stick r-f\},
\end{align*}
and $g_T=\hT_r+\stick r$ on $\{\Y_T>0\}$. Here $\hT_r$ is the passage time of $\BM$ from $-x$ to $r$, over the distance
$x+r$.
\end{corollary}

Both laws are passage times conditioned on later ones \eqref{eq:ptb}, and Proposition~\ref{prop:ptb} samples them. The
time $\hT_r-\hT_0$ is the operational time that the path spends away from $0$ between its first and last zeros. For $y=0$
the first line holds as written, and $g_T=T$.

\begin{proof}
Let $y>0$. Proposition~\ref{prop:zeros} with $s=T$ and $g=f+t+\stick r$, divided by
$\pYl(y,r,T\mid x)=2h(x+y+r,T-\stick r)$, gives the joint density
\[
	\frac{h(x,f)\,h(r,t)\,2h(y,\ T-\stick r-f-t)}{2h(x+y+r,\ T-\stick r)}
	~=~\frac{h(x,f)\,h(r+y,\ T-\stick r-f)}{h(x+r+y,\ T-\stick r)}\cdot\frac{h(r,t)\,h(y,\ T-\stick r-f-t)}{h(r+y,\ T-\stick r-f)},
\]
and the two factors are the two laws of the corollary. For $y=0$, the second line of Proposition~\ref{prop:zeros},
divided by $\stick h(x+r,T-\stick r)$, is the passage time conditioned on a later one with the levels $x$ and $r$.
\end{proof}

\subsection{Decomposition at the first and last zeros}\label{sec:bridge-given-zeros}

\begin{figure}[htb]
	\centering
	\includegraphics[width=\textwidth]{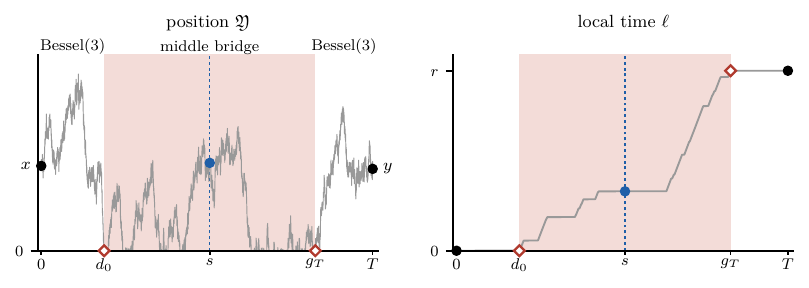}
    \caption[The bridge given position and local time]{The first and last zeros $d_0$ and $g_T$ delimit the middle bridge
    from $(0,0)$ to $(0,r)$ within the bridge from $(x,0)$ to $(y,r)$. On the left, the position follows a Bessel(3) bridge on each outer interval.
	The shaded middle bridge has operational length $g_T-d_0-\stick r$. On the right, the local time stays at
	$0$ up to $d_0$, increases to $r$ across the middle bridge, and stays at $r$ from $g_T$ to $T$.
	Algorithm~\ref{alg:bridge} draws the two zeros before sampling the value at $s$.}
	\label{fig:bridge}
\end{figure}

Given its two zeros, the path splits into three bridges (Figure~\ref{fig:bridge}). Before the first zero $d_0$ it is
a Bessel(3) bridge from $x$ to $0$ over $[0,d_0]$, with local time $0$. After the last zero $g_T$ it is a Bessel(3)
bridge from $0$ to $y$ over $[g_T,T]$, with local time $r$. Between the two zeros, the pair starts at $(0,0)$ at
time $d_0$, ends at $(0,r)$ at time $g_T$, and spends the operational time $\hT_r-\hT_0=g_T-d_0-\stick r$ away from
$0$. Its law depends on the zeros only through $r$ and this operational time. It is the only one of the three bridges
without an elementary law at one time, so we define it first, derive that law, and then identify it as the law of the
path between the zeros.

For $r,t>0$, the middle bridge of local time $r$ and operational length $t$ is the pair $(\Y,\locT)$ under
$\Probx{0}\{\iv\mid\hT_r=t\}$, on the physical times $[0,t+\stick r]$. It ends at $(0,r)$.
Appendix~\ref{app:bridges-paths} defines it as a law of paths, and Proposition~\ref{prop:bridge-path} shows that it is
also the bridge from $(0,0)$ to $(0,r)$ over the physical time $t+\stick r$. In the middle bridge, $s-\stick a$ and
$t-s+\stick a$ are the operational times spent away from $0$ before and after the time $s$, when the local time at $s$
is $a$.

\FloatBarrier

\begin{proposition}\label{prop:middle}
Let $r,t>0$, $0<s<t+\stick r$, $0<a<r$ and $z>0$. Then
\begin{align*}
	\Probx{0}\{\Y_s=0,\ \locT_s\in\md a\mid\hT_r=t\}
	~&=~\frac{\stick\,h(a,\ s-\stick a)\,h(r-a,\ t-s+\stick a)}{h(r,t)}\md a,\\
	\Probx{0}\{\Y_s\in\md z,\ \locT_s\in\md a\mid\hT_r=t\}
	~&=~\frac{2h(z+a,\ s-\stick a)\,h(z+r-a,\ t-s+\stick a)}{h(r,t)}\md z\md a .
\end{align*}
In particular, the law of $\locT_s$ under $\Probx{0}\{\iv\mid\hT_r=t\}$ has no atoms, and it is carried by the interval
$\big(\max\{0,(s-t)/\stick\},\ \min\{r,s/\stick\}\big)$, where both operational times $s-\stick a$ and $t-s+\stick a$ are
positive.
\end{proposition}

\begin{proof}
Let $\locT_s=a$. Then $\sigma_s=s-\stick a$, and by Lemma~\ref{lem:markov} at $s$ the time $\hT_r$ is $s-\stick a$ plus
the passage time of $\BM'$ over the distance $z+r-a$. With Proposition~\ref{prop:kernel} this gives
\[
	\Probx{0}\{\Y_s\in\md z,\ \locT_s\in\md a,\ \hT_r\in\md t\}~=~2h(z+a,\ s-\stick a)\,h(z+r-a,\ t-s+\stick a)\md z\md a\md t,
\]
and the same with $\stick h(a,\ s-\stick a)$ in place of $2h(z+a,\ s-\stick a)$ for $z=0$. Dividing by the density
$h(r,t)$ of $\hT_r$ gives the two displays.
\end{proof}

We now identify the middle bridge as the conditional law of the path between the first and last zeros, and the two
outer bridges as the Bessel(3) bridges of Section~\ref{sec:passage}. Let $f\geq0$, with $f=0$ if $x=0$, and $t>0$ with
$f+t+\stick r<T$, and put
$g\ass f+t+\stick r$. We write
\[
	\mathcal E~\ass~\{\hT_0=f,\ \hT_r-\hT_0=t,\ \Y_T=y,\ \locT_T=r\}
\]
for the event that fixes the zeros and the endpoint. On $\mathcal E$ the first zero is $f$ and, for $y>0$, the last zero
is $g$.

\begin{proposition}\label{prop:given-zeros}
Let $y>0$. If $s<f$, then $\locT_s=0$ and, for $z>0$,
\[
	\Probx{x}\{\Y_s\in\md z\mid\mathcal E\}~=~\frac{p_{\mathfrak A}(z,s\mid x)\,h(z,f-s)}{h(x,f)}\md z .
\]
If $s>g$, then $\locT_s=r$ and, for $z>0$,
\[
	\Probx{x}\{\Y_s\in\md z\mid\mathcal E\}~=~\frac{2h(z,s-g)\,p_{\mathfrak A}(y,T-s\mid z)}{2h(y,T-g)}\md z .
\]
If $f<s<g$, then $(\Y_s,\locT_s)$ given $\mathcal E$ has the law of the middle bridge of local time $r$ and operational
length $t$ at the time $s-f$.
\end{proposition}

The first two laws are Bessel(3) bridges of Lemma~\ref{lem:bessel}. The first goes from $x$ to $0$ over the time $f$, at
the time $s$. The second goes from $0$ to $y$ over the time $T-g$, at the time $s-g$. For $y=0$ the proposition holds
with $t=T-\stick r-f$. Then $g=T$, and the case $s>g$ does not occur.

\begin{proof}
We compute the joint density of $(\hT_0,\Y_s,\locT_s,g_T,\Y_T,\locT_T)$ under $\Probx{x}$ with the Markov property at $s$,
using Proposition~\ref{prop:zeros} on each side of $s$. Then we divide by the density $h(x,f)\,h(r,t)\,2h(y,T-g)$ of
$(\hT_0,g_T,\Y_T,\locT_T)$, which is Proposition~\ref{prop:zeros} at the time $T$.

Let $f<s<g$, $z>0$ and $0<a<r$. The first zero lies before $s$. Integrating Proposition~\ref{prop:zeros} at the time $s$
over $g_s$ by \eqref{eq:h-conv} gives the density $h(x,f)\,2h(z+a,\ s-f-\stick a)$ of $(\hT_0,\Y_s,\locT_s)$. After $s$,
the path is a new copy from $z$ whose local time increases by $r-a$, and its last zero before $T$ is $g_T$.
Proposition~\ref{prop:zeros} from $z$ over the time $T-s$, integrated over its first zero by \eqref{eq:h-conv}, gives
the density $h(z+r-a,\ g-s-\stick(r-a))\,2h(y,T-g)$ of $(g_T,\Y_T,\locT_T)$ given $(\Y_s,\locT_s)=(z,a)$. Since
$g-s-\stick(r-a)=t-(s-f)+\stick a$, the ratio is the second display of Proposition~\ref{prop:middle} at the time $s-f$.
The case $z=0$ is the same with $\stick h(a,\ s-f-\stick a)$ in place of $2h(z+a,\ s-f-\stick a)$, and it gives the
first display.

Let $s<f$. The state at $s$ is $(z,0)$ with density $p_{\mathfrak A}(z,s\mid x)$. Proposition~\ref{prop:zeros} from $z$
over the time $T-s$ gives the density $h(z,f-s)\,h(r,t)\,2h(y,T-g)$ of the rest, and the ratio is the first display.

Let $s>g$. Proposition~\ref{prop:zeros} at the time $s$ gives the density $h(x,f)\,h(r,t)\,2h(z,s-g)$ of
$(\hT_0,g_s,\Y_s,\locT_s)$. The last zero satisfies $g_T=g_s$ exactly when the path has no zero in $[s,T]$, and from $z$
this event has the density $p_{\mathfrak A}(y,T-s\mid z)$ at $\Y_T=y$. The ratio is the second display. For $y=0$, the
same computations hold with the atom $\stick h(r-a,\iv)$ of Proposition~\ref{prop:zeros} in place of
$h(z+r-a,\iv)\,2h(y,\iv)$, and with $g=T$.
\end{proof}

\FloatBarrier

\subsection{Sampling the middle bridge}\label{sec:middle}

Fix $r,t>0$ and a physical time $s\in(0,t+\stick r)$. We sample $(\Y_s,\locT_s)$ under $\Probx{0}\{\iv\mid\hT_r=t\}$.
As in Section~\ref{sec:forward}, we draw the local time first and then the position given the local time. In
Section~\ref{sec:forward} the local time was an explicit function of one half-normal. Here its distribution function is
explicit, and the draw inverts it. The local time at $s$ is at least $a$ exactly when it reached $a$ before $s$. This is
an event about the single passage time $\hT_a$, and the law of $\hT_a$ given $\hT_r=t$ is a passage time conditioned on
a later one (Section~\ref{sec:passage}).

\begin{proposition}\label{prop:middle-lt}
For $0<a<r$,
\[
	\Probx{0}\{\locT_s\geq a\mid\hT_r=t\}~=~\Probx{0}\{\hT_a\leq s-\stick a\mid\hT_r=t\}.
\]
As a function of $a$, the right side decreases continuously and strictly from $1$ at $a=\max\{0,(s-t)/\stick\}$ to $0$ at
$a=\min\{r,s/\stick\}$.
\end{proposition}

The right side is \eqref{eq:ptb-cdf} with the levels $a$ and $r-a$ and the time $t$, at $u=s-\stick a$. It is $0$ when
$s-\stick a\leq0$ and $1$ when $s-\stick a\geq t$. Hence, for a uniform $\uni$, the root $a$ of
\[
	\Probx{0}\{\hT_a\leq s-\stick a\mid\hT_r=t\}~=~\uni
\]
has the law of $\locT_s$ under $\Probx{0}\{\iv\mid\hT_r=t\}$.

\begin{proof}
The local time is nondecreasing and first reaches $a$ at the physical time $\hT_a+\stick a$. Hence
$\{\locT_s\geq a\}=\{\hT_a+\stick a\leq s\}$ for every path, which is the identity. It remains to check that the two
sides agree for the versions we use, those of Proposition~\ref{prop:middle} on the left and of \eqref{eq:ptb} on the
right, at every admissible $t$. By the strong Markov property of $\BM$ at
$\hT_a$, the times $\hT_a$ and $\hT_r-\hT_a$ are independent with densities $h(a,\iv)$ and $h(r-a,\iv)$. This gives
\[
	\Probx{0}\{\hT_a\leq s-\stick a,\ \hT_r\in\md t\}~=~\int_0^{s-\stick a}h(a,w)\,h(r-a,t-w)\md w\md t .
\]
For $t>s-\stick r$ the local time at $s$ is below $r$, and the joint density in the proof of
Proposition~\ref{prop:middle} gives
\begin{align*}
	\Probx{0}\{\locT_s\geq a,\ \hT_r\in\md t\}
	~&=~\int_a^r\Big(\stick h(b,s-\stick b)\,h(r-b,t-s+\stick b)\\
	&\qquad\qquad+\int_0^\infty2h(z+b,s-\stick b)\,h(z+r-b,t-s+\stick b)\md z\Big)\md b\md t .
\end{align*}
The two right sides are densities in $t$ of the same measure. Both are continuous in $t$, so they agree for every $t$,
and dividing by $h(r,t)$ gives the identity. Its right side, evaluated by \eqref{eq:ptb-cdf}, is continuous in $a$. By
the last display, its derivative in $a$ is minus the total mass of the two displays of Proposition~\ref{prop:middle} at
the level $a$, which is negative.
\end{proof}

Given the local time, the position has an atom at $0$ and a density on $(0,\infty)$. The atom and the density are
known up to a common normalizing constant, which we do not evaluate. We sample the two together by rejection, from a
proposal that keeps the atom as it is and dominates the density by Gamma densities. The operational times
$s-\stick a$ and $t-s+\stick a$ spent away from $0$ before and after $s$ are both positive. We put
\[
	v~\ass~\frac{(s-\stick a)(t-s+\stick a)}{t},\qquad
	k~\ass~\sqrt v\,\Big(\frac a{s-\stick a}+\frac{r-a}{t-s+\stick a}\Big).
\]
Here $v$ is the variance at the time $s-\stick a$ of a Brownian bridge over the time $t$, and $k>0$.

\begin{proposition}\label{prop:middle-pos}
For $0<a<r$ with $0<s-\stick a<t$ and $z>0$,
\begin{align*}
	\Probx{0}\{\Y_s\in\md z\mid\locT_s=a,\ \hT_r=t\}
	~&=~\Probx{0}\{\Y_s=0\mid\locT_s=a,\ \hT_r=t\}\\
	&\qquad\cdot\frac{2}{\stick a(r-a)}\,(z+a)(z+r-a)\,\exp\Big\{-\frac{kz}{\sqrt v}-\frac{z^2}{2v}\Big\}\md z .
\end{align*}
\end{proposition}

In the variable $w=z/\sqrt v$, the atom and the density are proportional, with the same constant, to
\[
	\tfrac12\stick a(r-a)\qquad\text{and}\qquad\sqrt v\,\big(a(r-a)+r\sqrt v\,w+vw^2\big)\me^{-kw-w^2/2}.
\]

\begin{proof}
The conditional atom and density given $\locT_s=a$ are the two displays of Proposition~\ref{prop:middle} divided by the
same density of $\locT_s$. Their ratio is therefore the ratio of the displays. Since
$h(c+z,u)/h(c,u)=\frac{c+z}c\exp\{-(2cz+z^2)/2u\}$ and $\frac1{s-\stick a}+\frac1{t-s+\stick a}=\frac1v$,
\[
	\frac{2h(z+a,\ s-\stick a)\,h(z+r-a,\ t-s+\stick a)}{\stick h(a,\ s-\stick a)\,h(r-a,\ t-s+\stick a)}
	~=~\frac2\stick\cdot\frac{z+a}a\cdot\frac{z+r-a}{r-a}\cdot\exp\Big\{-\frac{kz}{\sqrt v}-\frac{z^2}{2v}\Big\}.
\]
The form in $w$ follows from $(z+a)(z+r-a)=a(r-a)+rz+z^2$ and $\md z=\sqrt v\md w$.
\end{proof}

In the variable $w$, the density on $(0,\infty)$ is a combination of the three functions $w^j\me^{-kw-w^2/2}$,
$j=0,1,2$, with positive coefficients. A multiple of a Gamma density dominates each of the three, and we sample the
Gamma law as a sum of exponentials. With the rate and the mass
\[
	\lambda_j~\ass~\tfrac12\big(k+\sqrt{k^2+4(j+1)}\big),\qquad I_j~\ass~\int_0^\infty w^j\me^{-kw-w^2/2}\md w,
\]
Lemma~\ref{lem:envelope} gives
\[
	w^j\me^{-kw-w^2/2}~\leq~\me^{(\lambda_j-k)^2/2}\,j!\,\lambda_j^{-j-1}\cdot\frac{\lambda_j^{j+1}w^j\me^{-\lambda_jw}}{j!},\qquad w>0,
\]
where the last factor is the density of $\GammaLaw(j+1,\lambda_j)$. A proposal $w$ from this Gamma law is accepted
with the probability $\exp\{-(w-\lambda_j+k)^2/2\}$, the ratio of the two sides, and the mean acceptance probability
of the proposal is
\[
	A_j(k)~\ass~\frac{\lambda_j^{j+1}\me^{-(\lambda_j-k)^2/2}}{j!}\,I_j~\geq~A_2(0)~=~0.72656\ldots,\qquad k\geq0,\ j\leq2,
\]
uniformly in the parameters.

\subsection{The bridge algorithm}\label{sec:bridge-alg}

Algorithm~\ref{alg:bridge} combines the previous subsections. It draws the zeros $f$ and $g$ by
Corollary~\ref{cor:zeros}, and it returns them with the value at $s$. If $s<f$ or $s>g$, the value is the Bessel(3)
bridge of Proposition~\ref{prop:given-zeros} at $s$, by Lemma~\ref{lem:bessel}. Otherwise $s$ falls in the middle
bridge, at its time $s'\ass s-f$. The local time is drawn by inversion (Proposition~\ref{prop:middle-lt}), and a loop
then draws the position from Proposition~\ref{prop:middle-pos}, with the envelopes of Lemma~\ref{lem:envelope}. With
$s'$ in place of $s$ in $v$ and $k$, the loop uses
\begin{gather*}
	c_0\ass a(r-a),\qquad c_1\ass r\sqrt v,\qquad c_2\ass v,\\
	m_j\ass\sqrt v\,c_j\,\me^{(\lambda_j-k)^2/2}\,j!\,\lambda_j^{-j-1},\qquad M\ass\tfrac12\stick a(r-a)+m_0+m_1+m_2 .
\end{gather*}

\begin{algorithm}[H]
\caption{The bridge at one time (Theorem~\ref{thm:bridge})}\label{alg:bridge}
\begin{algorithmic}[1]
\Require ends $x,y\geq0$, a local time $r>0$, a length $T>\stick r$ and a time $s\in(0,T)$
\State draw the first zero $f$ and the last zero $g$ by Corollary~\ref{cor:zeros} and Proposition~\ref{prop:ptb}
\State \textbf{if} $s<f$ \Return $(f,g,\bessb_s,0)$, from \eqref{eq:bessel-sampler} with $(z,z',t,u)=(x,0,f,s)$
\State \textbf{if} $s>g$ \Return $(f,g,\bessb_{s-g},r)$, from \eqref{eq:bessel-sampler} with $(z,z',t,u)=(0,y,T-g,s-g)$
\State put $s'\ass s-f$ and $t\ass g-f-\stick r$, draw $\uni$, and let $a$ be the root of
$\Probx{0}\{\hT_a\leq s'-\stick a\mid\hT_r=t\}=\uni$\Comment{\eqref{eq:ptb-cdf}}
\State compute $v$, $k$, $\lambda_j$, $c_j$, $m_j$ for $j\leq2$, and $M$, from $a$, $r$, $t$ and $s'$
\Loop
	\State draw $\uni$ and put $b\ass\uni M$; \textbf{if} $b<\frac12\stick a(r-a)$ \Return $(f,g,0,a)$
	\State let $j$ be the least index with $b<\frac12\stick a(r-a)+m_0+\dots+m_j$
	\State draw $\erv_1,\dots,\erv_{j+1}$ and $\uni$, and put $w\ass(\erv_1+\dots+\erv_{j+1})/\lambda_j$
	\State \textbf{if} $\uni\leq\exp\{-(w-\lambda_j+k)^2/2\}$ \Return $(f,g,\sqrt v\,w,\ a)$
\EndLoop
\end{algorithmic}
\end{algorithm}

The theorem states the joint law of the value and the zeros, which is what the refinement of
Section~\ref{sec:refine} uses, and the law \eqref{eq:bridge-one} of the value alone follows from it. The root $a$ in
the theorem is the exact root of the equation of Proposition~\ref{prop:middle-lt}; the implementation computes it by
bisection (Section~\ref{sec:numerics}).

\begin{theorem}\label{thm:bridge}
Let $T>0$, $x,y\geq0$ and $r>0$ with $\stick r<T$, and let $s\in(0,T)$. The zeros and the value $(f,g,Y,L)$ returned by
Algorithm~\ref{alg:bridge} satisfy
\[
	(f,\,g,\,Y,\,L)~\law~(\hT_0,\,g_T,\,\Y_s,\,\locT_s)\qquad\text{under }\Probx{x}\{\iv\mid\Y_T=y,\ \locT_T=r\}.
\]
In particular, $(Y,L)$ has the law \eqref{eq:bridge-one},
\begin{align*}
	\Prob\{Y\in\md z,\ L\in\md a\}
	~&=~\frac{\pYl(z,a,s\mid x)\,\pYl(y,r-a,T-s\mid z)}{\pYl(y,r,T\mid x)}\\
	&\qquad\cdot\big(\delta_0(\md z)+\md z\big)\big(\delta_0(\md a)+\md a+\delta_r(\md a)\big).
\end{align*}
\end{theorem}

\begin{proof}
By Corollary~\ref{cor:zeros}, $(f,g)$ has the law of $(\hT_0,g_T)$ under the bridge. Given the zeros, $(Y,L)$ has the
conditional law of $(\Y_s,\locT_s)$ given the zeros and the endpoint. For $s<f$ and $s>g$, this follows from
Proposition~\ref{prop:given-zeros} and Lemma~\ref{lem:bessel}. For $f<s<g$, it follows from
Proposition~\ref{prop:given-zeros}, with $L$ drawn by Proposition~\ref{prop:middle-lt} and $Y$ by the loop. In one pass,
the loop returns $Y=0$ with probability $\frac12\stick a(r-a)/M$. By Lemma~\ref{lem:envelope}, it chooses $j$, proposes
$w$ in $\md w$ and accepts with probability
\[
	\frac{m_j}M\cdot\frac{\lambda_j^{j+1}w^j\me^{-\lambda_jw}}{j!}\cdot\me^{-(w-\lambda_j+k)^2/2}\md w
	~=~\frac1M\,\sqrt v\,c_j\,w^j\me^{-kw-w^2/2}\md w .
\]
Summed over $j$, a pass returns a value in $\sqrt v\md w$ with probability
$\frac1M\sqrt v\,(c_0+c_1w+c_2w^2)\me^{-kw-w^2/2}\md w$. These are the weights of Proposition~\ref{prop:middle-pos} with
the common constant $1/M$. The passes are independent, so the value returned has the conditional law of $\Y_s$ given
$\locT_s=a$ \citep{vonNeumann51}. Combining the law of the zeros with the law of the value given the zeros gives the
first display, and integrating out the zeros gives the second.
\end{proof}

The number of passes of the loop in Algorithm~\ref{alg:bridge} has a geometric bound that does not depend on the parameters. We write
$K$ for this number when $s$ falls in the middle bridge.

\begin{proposition}\label{prop:cost}
Algorithm~\ref{alg:bridge} uses finitely many elementary random variables almost surely. Each pass of the loop returns a value with
probability at least $A_2(0)=0.72656\ldots$, whatever $x$, $y$, $r$, $T$, $s$ and the drawn $a$. Hence
\[
	\Prob\{K>n\}~\leq~\big(1-A_2(0)\big)^n,\qquad\E[K]~\leq~\frac1{A_2(0)}~=~1.3763\ldots
\]
\end{proposition}

\begin{proof}
A pass succeeds with probability
\[
	\frac{\frac12\stick a(r-a)+\sum_jc_j\sqrt v\,I_j}{\frac12\stick a(r-a)+\sum_jm_j}~\geq~\min_{j\leq2}\frac{c_j\sqrt v\,I_j}{m_j}
	~=~\min_{j\leq2}A_j(k)~\geq~A_2(0),
\]
by Lemma~\ref{lem:envelope}, since $k\geq0$. The passes are independent given $a$, so $K$ is geometric with success
probability at least $A_2(0)$. Besides the passes, the zeros use at most two normals and four uniforms. A Bessel(3) value
uses one normal and one exponential, and one uniform when both of its ends are positive. The local time uses one
uniform, and each pass uses two uniforms and at most three exponentials.
\end{proof}

Since $k>0$, the bound $A_2(0)$ is not attained. It is approached as $r/\sqrt t\to0$.

When $r=0$ and $x,y>0$, the path has no zero. Its local-time gains on
$[0,s]$ and $[s,T]$ are therefore both zero. Applying the Markov
property at $s$ and the zero-gain case of Proposition~\ref{prop:kernel}
on each interval gives
\[
    \Probx{x}\{\Y_s\in\md z\mid\Y_T=y,\ \locT_T=0\}
    ~=~
    \frac{p_{\mathfrak A}(z,s\mid x)\,p_{\mathfrak A}(y,T-s\mid z)}
         {p_{\mathfrak A}(y,T\mid x)}\,\md z,
    \qquad z>0.
\]
This is the Bessel(3) bridge from $x$ to $y$ of
Lemma~\ref{lem:bessel}, sampled directly by \eqref{eq:bessel-sampler}.
If $r=0$ and $xy=0$, or if $\stick r\geq T$, then
$\pYl(y,r,T\mid x)=0$ and \eqref{eq:bridge-def} does not define the
bridge.

\begin{remark}
As $\stick\to0$, the middle bridge spends no time at $0$, and the atom in Proposition~\ref{prop:middle-pos}
disappears. Proposition~\ref{prop:given-zeros} then becomes the law of reflected Brownian motion and its local time
inside a bridge given both. For $f<s<g$, the second display of Proposition~\ref{prop:middle} then gives the joint density of the reflected position and its local time inside the middle bridge.
\citet{papaspiliopoulosExactSamplingDiffusions2016} give the corresponding density for signed Brownian motion as a product of two Rayleigh factors after a change of variables.
Taking the absolute value of their position variable gives our formula.
Stickiness replaces the physical times in the arguments by
the operational times $s-\stick a$ and $t-s+\stick a$, and it adds the atom at $0$.
\end{remark}

\section{Refinement in any order}\label{sec:refine}

We refine a sampled path by querying it at a new time. The new value must have its conditional law given everything
sampled so far, including the zeros that earlier queries drew. 
Suppose a query inside the middle bridge finds the
path above zero. The query also draws the two zeros bounding the excursion through that time. If the next query falls 
between the first of these zeros and the sampled point, we know
that the path stays positive between them. We therefore draw from the Bessel(3) bridge from zero to the sampled
position, keeping the local time at its sampled level. Figure~\ref{fig:query} explains this in one picture, for a single path. 

Following the terminology of exact simulation \citep{beskosExactSimulationDiffusions2005,pollockExactEstrongSimulation2016},
we call the sampled values and the recorded zeros a skeleton. The next query time may itself be chosen from the
skeleton. To justify successive draws, we describe the conditional law of the path given the skeleton and show how
that law changes when a query adds new points.

\subsection{Skeletons and queries}\label{sec:records}

We list the points of a skeleton $\rec$ in increasing order of time as $(s_i,z_i,a_i)$, $0\leq i\leq m$.
Here $z_i\geq0$ is the position at $s_i$ and $a_i\geq0$ is the local time. For consecutive points put
\[
    \Delta_i\ass s_{i+1}-s_i,\qquad \delta_i\ass a_{i+1}-a_i,\qquad 0\leq i<m.
\]
The local times are nondecreasing. We call $(s_i,s_{i+1})$ a constant interval if $\delta_i=0$ and an increasing
interval if $\delta_i>0$. In an increasing interval both endpoint positions are zero and $\Delta_i>\stick\delta_i$.
Of the physical duration $\Delta_i$, the path spends $\stick\delta_i$ at zero and the remaining
$\Delta_i-\stick\delta_i$ away from zero. Appendix~\ref{app:kernel} gives the full conditions for admissibility.

For a finite set $\mathcal T$ of times containing $0$, let $\rrec(\mathcal T)$ denote the skeleton read from
the path at those times. Recall the two times from Section~\ref{sec:markov},
\[
    d_v\ass\inf\{t\geq v:\Y_t=0\},\qquad
    g_v\ass\sup\{0\leq t\leq v:\Y_t=0\},
\]
with $d_v=\infty$ and $g_v=0$ for the respective empty sets. These are the first zero at or after $v$
and the last zero at or before $v$ whenever those zeros exist. In particular, $g_v=d_v=v$ if $\Y_v=0$.

To form the skeleton, start with the observations and add zeros only between consecutive times
$u<u'$ of $\mathcal T$ for which $\locT_u<\locT_{u'}$. With $(u,u')$ ranging over these pairs, put
\begin{equation}\label{eq:revealed}
    \begin{aligned}
    \rrec(\mathcal T)
      &\ass \{(v,\Y_v,\locT_v):v\in\mathcal T\}\\
       &{}\cup\{(d_u,0,\locT_u):\Y_u>0\}
          \cup\{(g_{u'},0,\locT_{u'}):\Y_{u'}>0\}.
    \end{aligned}
\end{equation}
The additional points are the recorded zeros. To see why they carry the indicated local times, consider a pair
$(u,u')$ across which the local time increases. It stays at $\locT_u$ from $u$ to $d_u$ and at $\locT_{u'}$ from
$g_{u'}$ to $u'$. All the gain therefore occurs between $d_u$ and $g_{u'}$. An endpoint already at zero
contributes no additional point. Almost surely $\rrec(\mathcal T)$ is a skeleton and its constant intervals
have no zero in their interiors (Appendix~\ref{app:revealed}). A later query may select a recorded zero.
Its time then belongs to the enlarged set of query times.

When only the values at the times of $\mathcal T$ are given, we draw the additional zeros in
\eqref{eq:revealed} conditionally on those values. We call this operation completion. Consider consecutive
observations $(u_1,z_1,a_1)$ and $(u_2,z_2,a_2)$ with $\delta\ass a_2-a_1>0$ and $z_1+z_2>0$.
Write $d$ for the first zero and $g$ for the last zero in $[u_1,u_2]$. Set $d=u_1$ if $z_1=0$ and
$g=u_2$ if $z_2=0$. By Corollary~\ref{cor:zeros}, if $z_1>0$, we draw $d-u_1$ from
\[
    \Probx{0}\{\hT_{z_1}\in\iv\mid\hT_{z_1+\delta+z_2}=u_2-u_1-\stick\delta\}.
\]
Once $d$ is known, if $z_2>0$, we draw $g-d-\stick\delta$ from
\[
    \Probx{0}\{\hT_\delta\in\iv\mid\hT_{\delta+z_2}=u_2-d-\stick\delta\}.
\]
We add $(d,0,a_1)$ if $z_1>0$ and $(g,0,a_2)$ if $z_2>0$. Given the observations, the draws for different
pairs are independent.

A query at $u$ returns a skeleton $\qry_u(\rec)$, a measurable function of $(\rec,u)$ and fresh
elementary random variables independent of all previous draws. At a time already recorded it returns $\rec$.
Beyond the last time $s_m$ of an unconditioned path, Algorithm~\ref{alg:forward} draws the position $z$ and
local-time gain $a$ from $z_m$ over $u-s_m$. We add $(u,z,a_m+a)$ and complete the new pair.

Suppose now that $s_i<u<s_{i+1}$. In a constant interval we draw the position at $u-s_i$ from the Bessel(3)
bridge of Lemma~\ref{lem:bessel} from $z_i$ to $z_{i+1}$ of duration $\Delta_i$. Its local time remains $a_i$.
In an increasing interval we draw $(z,a)$ at $u-s_i$ from the middle bridge of local-time gain $\delta_i$ and
operational length $\Delta_i-\stick\delta_i$, using Section~\ref{sec:middle}. We add $(u,z,a_i+a)$.
If $z=0$, this single point splits the interval into two middle bridges with gains $a$ and $\delta_i-a$.

\begin{figure}[H]
	\centering
	\includegraphics[width=\textwidth]{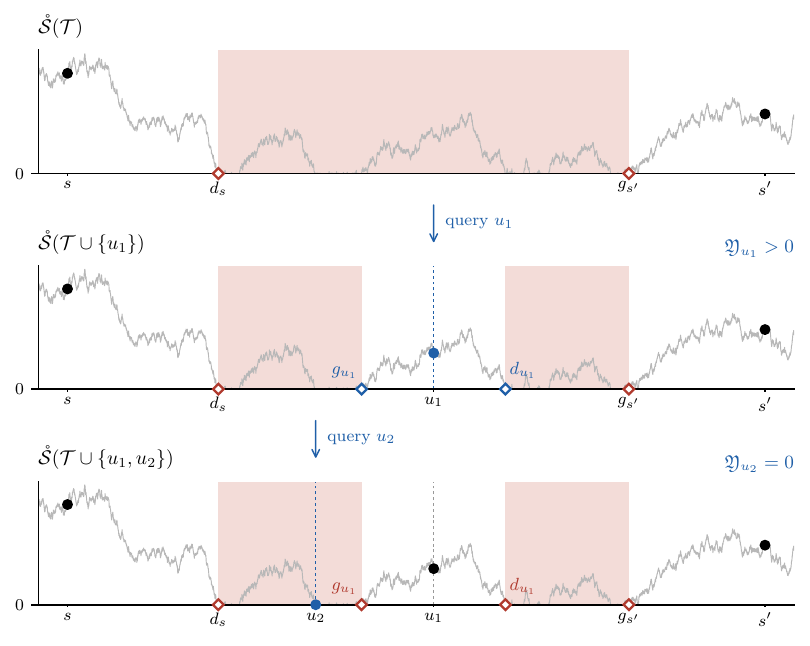}
	\caption[Two queries on one path]{The top panel shows the skeleton between consecutive query times $s<s'$. Its
recorded zeros $d_s$ and $g_{s'}$ bound the shaded middle bridge. In the middle
panel the path is positive at $u_1$, so the query adds its value and the two
zeros bounding its excursion. In the bottom panel the path is zero at $u_2$,
so adding that point splits the remaining middle bridge containing $u_2$ in
two. In each step, with $\mathcal T$ the set of query times before the new
query at $u$,
\[
    \rrec(\mathcal T\cup\{u\})~=~\rrec(\mathcal T)\cup\big\{(g_u,0,\locT_u),\,(u,\Y_u,\locT_u),\,(d_u,0,\locT_u)\big\},
\]
where $g_u=d_u=u$ when $\Y_u=0$.}
	\label{fig:query}
\end{figure}

If $z>0$, let $g$ and $d$ be the zeros bounding the excursion through $u$. Almost surely $0<a<\delta_i$
(Proposition~\ref{prop:middle}). Given $(z,a)$ and the old skeleton, the backward duration $u-g$ and the
forward duration $d-u$ are independent. Their respective conditional laws are
\begin{gather*}
    \Probx{0}\{\hT_z\in\iv\mid\hT_{z+a}=u-s_i-\stick a\},\\
    \Probx{0}\{\hT_z\in\iv\mid\hT_{z+\delta_i-a}=\Delta_i-\stick\delta_i-(u-s_i-\stick a)\}.
\end{gather*}
Both zeros have local time $a_i+a$. We add $(g,0,a_i+a)$ and $(d,0,a_i+a)$ to the skeleton
(Figure~\ref{fig:query}).

\subsection{The conditional law given a skeleton}\label{sec:record-law}

The next query might be the midpoint of an interval whose ends were revealed by earlier queries. For each
fixed $u$, a conditional identity may fail on a set of skeletons of probability zero. Since this set can
depend on $u$, a time chosen from the skeleton can select precisely such an exception. We therefore construct
the conditional law of the whole path and prove the query identity for every admissible pair $(\rec,u)$.

Write $\Bes_{z,z',t}$ for the law of a Bessel(3) bridge from $z$ to $z'$ of duration $t$ paired with zero local-time increment.
Write $\Mid_{\delta,t}$ for the joint law of a middle bridge with local-time gain $\delta$ and operational duration $t$.
Appendix~\ref{app:bridges-paths} constructs both laws.

Let $\Ker_\rec$ denote the conditional path law given the skeleton $\rec$. Under this law the bridges on successive
intervals are independent. On $[s_i,s_{i+1}]$ the path $(\Y_{s_i+v},\locT_{s_i+v}-a_i)_{0\leq v\leq\Delta_i}$ has law
\[
    \begin{cases}
        \Bes_{z_i,z_{i+1},\Delta_i}, & \delta_i=0,\\
        \Mid_{\delta_i,\Delta_i-\stick\delta_i}, & \delta_i>0.
    \end{cases}
\]
Each bridge starts at time $s_i$ with local time $a_i$. For an unconditioned process, append an independent
copy from $z_m$ after $s_m$, adding $a_m$ to its local time. For a bridge on $[0,T]$, the last point is at
$T$ and the construction ends there. Appendix~\ref{app:kernel} defines these kernels on the admissible
skeletons $\Recs$ and $\Recs_T$, respectively. The skeletons read from the path belong to the corresponding
set almost surely. In particular, $\Ker_{\{(0,x,0)\}}=\Probx{x}$.

Let $\mathcal T_\rec$ be the set of times of $\rec$. Given a path $w$ and a finite set of additional
times $\mathcal V$, write
\[
    \rec_{\mathcal V}(w)\ass\rrec(\mathcal T_\rec\cup\mathcal V)(w),
    \qquad \rec_u(w)\ass\rec_{\{u\}}(w).
\]
Thus $\rec_u(w)$ is read from a path, whereas $\qry_u(\rec)$ is drawn by the sampler. Write
$\qry_u(\rec,\cdot)$ for the law of that draw. It is a kernel in $(\rec,u)$ and is the point mass at $\rec$
when $u\in\mathcal T_\rec$. Proposition~\ref{prop:record} gives the density of $\rec_{\mathcal V}$ under
$\Ker_\rec$ in terms of passage-time densities $h$ and killed transition densities $p_{\mathfrak A}$.

\begin{proposition}\label{prop:path-query}
Let $\rec\in\Recs$ and $u>0$ in the unconditioned case, or
$\rec\in\Recs_T$ and $0<u<T$ in the bridge case. For all bounded measurable
functions $F$ of the path and $G$ of a skeleton,
\[
	\int F(w)\,G\big(\rec_u(w)\big)\,\Ker_\rec(\md w)~=~\int\Ker_{\rec'}[F]\;G(\rec')\;\qry_u(\rec,\md\rec').
\]
\end{proposition}

Taking $F=1$ shows that a query at $u$ draws exactly what a path of law
$\Ker_\rec$ would reveal there. The full identity gives the conditional law of that path after the query.
If $\rec'$ is the skeleton revealed, $\Ker_{\rec'}$ is a version of this conditional law.
The proof in Appendix~\ref{app:refine-proof} uses the factorization of Proposition~\ref{prop:record}.

\subsection{Repeated queries}\label{sec:repeated}

We make queries one after another, each at a time chosen from the skeleton so far. Under $\Probx{x}$ the input consists of the values of $(\Y,\locT)$ at a fixed finite set of times
$\mathcal T_0$ containing $0$. These values may be sampled by Theorem~\ref{thm:forward}. Under the bridge the
input consists of the endpoints $(0,x,0)$ and $(T,y,r)$, with $r\geq0$ and $\pYl(y,r,T\mid x)>0$. Here
$\mathcal T_0=\{0,T\}$ and we use the path law of Proposition~\ref{prop:bridge-path}, whose finite-dimensional
laws are those of Section~\ref{sec:bridge}. When $r=0$, the positivity condition requires $x,y>0$ and the
input is a single constant interval. In either case let $\rec_0$ denote the completed input.

At step $j$ a measurable function $\chi_j$ chooses the time from the current skeleton and an auxiliary uniform
variable $\mathfrak v_j$. The auxiliary variables are independent of one another, of the path, and of the
random variables used for sampling. The function $\chi_j$ takes values in $(0,\infty)$ for the unconditioned process and
in $(0,T)$ for the bridge. Thus Algorithm~\ref{alg:refine} has the recursion
\begin{equation}\label{eq:alg-iter}
    u_j\ass\chi_j(\rec_{j-1},\mathfrak v_j),\qquad
    \rec_j\ass\qry_{u_j}(\rec_{j-1}),\qquad j\geq1.
\end{equation}
Write $(Y_j,L_j)$ for the position and local time returned at $u_j$.

For comparison, apply the same functions $\chi_j$ to the successive skeletons of one path. At step $j$ this gives
\begin{equation}\label{eq:path-iter}
    u'_j\ass\chi_j\big(\rrec(\mathcal T_{j-1}),\mathfrak v_j\big),\qquad
    \mathcal T_j\ass\mathcal T_{j-1}\cup\{u'_j\},\qquad j\geq1.
\end{equation}
Let $\mathcal G_k$ be the $\sigma$-field generated by $\rrec(\mathcal T_0),\ldots,\rrec(\mathcal T_k)$
and $\mathfrak v_1,\ldots,\mathfrak v_k$, and let $\E$ denote expectation under $\Probx{x}$ or under the bridge.

\begin{algorithm}[H]
\caption{Refinement in any order (Theorem~\ref{thm:refine})}\label{alg:refine}
\begin{algorithmic}[1]
\Require the completed input $\rec_0$ and measurable functions $\chi_1,\chi_2,\dots$
\For{$j=1,2,\dots$}
    \State draw $\mathfrak v_j$ uniform on $(0,1)$, put $u\ass\chi_j(\rec_{j-1},\mathfrak v_j)$ and $\rec\ass\rec_{j-1}$
    \If{$u$ is the time of a point of $\rec$}
        \State leave $\rec$ unchanged
    \ElsIf{$u$ is after the last point $(s_m,z_m,a_m)$ of $\rec$}
        \State draw $(z,a)$ by Algorithm~\ref{alg:forward} from $z_m$ over the time $u-s_m$, add $(u,z,a_m+a)$, and
        complete the new pair
    \ElsIf{$u$ lies in a constant interval $(s_i,s_{i+1})$ of $\rec$}
        \State add $(u,\bessb_{u-s_i},a_i)$, from \eqref{eq:bessel-sampler} with $(z,z',t)=(z_i,z_{i+1},\Delta_i)$
    \Else\Comment{an increasing interval $(s_i,s_{i+1})$}
        \State draw $(z,a)$ from the middle bridge of local time $\delta_i$ and operational length
        $\Delta_i-\stick\delta_i$ at its time $u-s_i$, by the local time and the loop of Section~\ref{sec:middle}, and add
        $(u,z,a_i+a)$
        \State \textbf{if} $z>0$, draw and add the two recorded zeros of Section~\ref{sec:records}
    \EndIf
    \State $\rec_j\ass\rec$
\EndFor
\end{algorithmic}
\end{algorithm}

\begin{theorem}\label{thm:refine}
Under $\Probx{x}$, or under the bridge, for every $k$ and every bounded measurable function $F$ of the path,
\begin{align*}
    \big(\rec_0,\rec_1,\dots,\rec_k\big)~&\law~\big(\rrec(\mathcal T_0),\rrec(\mathcal T_1),\dots,\rrec(\mathcal T_k)\big),\\
    \big(u_j,\,Y_j,\,L_j\big)_{j\leq k}~&\law~\big(u'_j,\,\Y_{u'_j},\,\locT_{u'_j}\big)_{j\leq k},\\
    \E\big[F(\Y,\locT)\bigm|\mathcal G_k\big]~&=~\Ker_{\rrec(\mathcal T_k)}[F].
\end{align*}
\end{theorem}

\begin{figure}[htb]
    \centering
    \includegraphics[width=\textwidth]{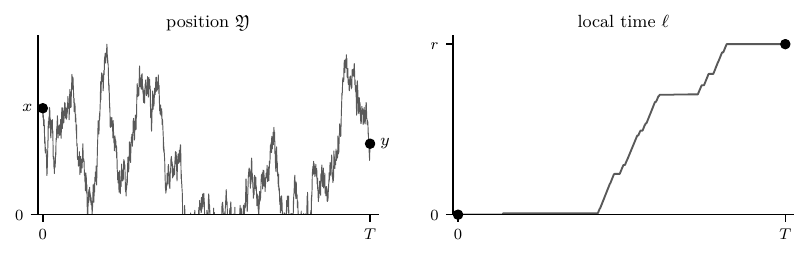}
    \caption[A refined bridge]{One sample of the bridge of $(\Y,\locT)$ from $(0.6,0)$ to $(0.4,0.8)$ over $[0,2]$,
    with $\stick=0.5$, after $3000$ queries by Algorithm~\ref{alg:refine} at the midpoint of the longest interval (seed 20261116). The
    skeleton has $3096$ points, which include the recorded zeros, and its values come from the implementation of
    Section~\ref{sec:numerics}. The lines join them for display. The local time increases only on the zero set.}
    \label{fig:refined}
\end{figure}

The first two identities give the joint laws of the successive skeletons and of the times and values
returned. The last identifies the current skeleton as sufficient to determine the conditional law of the path
given the entire query history. Appendix~\ref{app:refine-proof} proves the theorem by induction on $k$, using
Proposition~\ref{prop:path-query} at each step. Figure~\ref{fig:refined} shows one choice of the functions $\chi_j$. 
Each query bisects the longest interval of the current skeleton, and the recorded zeros are among the ends of those intervals.

Constant functions $\chi_j$ give refinement at prescribed times.

\begin{corollary}\label{cor:moves}
Let $\mathcal T$ be finite with $0\in\mathcal T$, and let $u>0$ with $u\notin\mathcal T$. Under $\Probx{x}$, and under the bridge if
$\{0,T\}\subseteq\mathcal T\subseteq[0,T]$ and $u<T$,
\[
    \qry_u\big(\rrec(\mathcal T)\big)~\law~\rrec\big(\mathcal T\cup\{u\}\big).
\]
\end{corollary}

\begin{proof}
Under $\Probx{x}$, take $\mathcal T_0=\mathcal T$ and $\chi_1\equiv u$ in Theorem~\ref{thm:refine}. Then $\rec_0\law\rrec(\mathcal T)$ and
$\qry_u(\rec_0)=\rec_1\law\rrec(\mathcal T\cup\{u\})$. The query uses fresh variables, so the law of $\qry_u(\rec_0)$ depends on
$\rec_0$ only through its law, and the identity follows. Under the bridge, queries at the times of $\mathcal T\setminus\{0,T\}$ in
increasing order give a skeleton with the law of $\rrec(\mathcal T)$, and the same argument applies.
\end{proof}

For a forward query beyond the last observation, the corollary recovers Theorem~\ref{thm:forward}.
From the two bridge endpoints it recovers Theorem~\ref{thm:bridge} when $r>0$, with $d_0$ and $g_T$
the first and last zeros used by Algorithm~\ref{alg:bridge}.

The skeleton also records how much time the path spent at zero. On $(s_i,s_{i+1})$ that time is
$\stick\delta_i$, and almost surely the open interval contains a zero exactly when $\delta_i>0$.
Zeros at the endpoints are read directly from $z_i$ and $z_{i+1}$.

\begin{remark}
Completion may be deferred until the first query between a pair of input points. Until then the conditional law on the interval
is the endpoint bridge of Section~\ref{sec:bridge}. Proposition~\ref{prop:bridge-path} expresses its
law as the mixture over its first and last zeros. Drawing the required zeros at that first query therefore gives the same conditional law for its answer.
The comparison uses the same function $\chi_j$ in both constructions, based only on the information available
before the deferred zeros are drawn.
\end{remark}

\FloatBarrier

\section{Implementation and a numerical check}\label{sec:numerics}

\begin{figure}[htb]
	\centering
	\includegraphics[width=\textwidth]{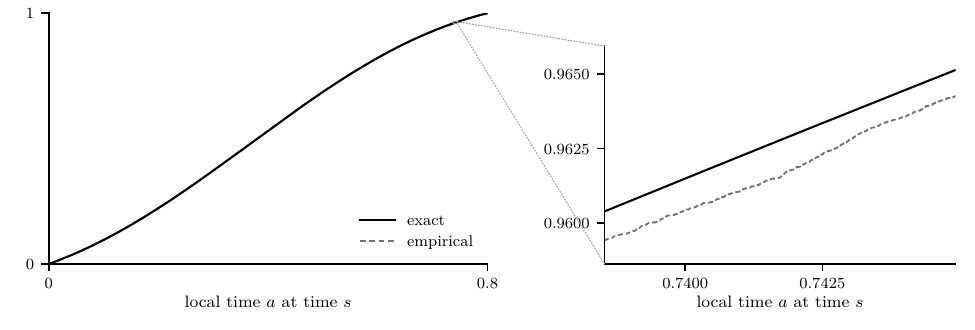}
	\caption[The local time of a middle bridge]{The local time at $s=0.9$ for the bridge from $(0,0)$ to $(0,0.8)$ over $[0,2]$
	with $\stick=0.5$. The left panel compares its exact distribution function
	(Proposition~\ref{prop:middle-lt}) with the empirical distribution function from $10^5$ independent
	runs of Algorithm~\ref{alg:bridge} (seed 20261117). The right panel enlarges the boxed region around
	$a=0.742$, where the largest observed difference is $0.0013$.}
	\label{fig:lt-check}
\end{figure}

The code, the scripts of Figures~\ref{fig:refined} and~\ref{fig:lt-check} and two interactive pages are available
at \url{https://github.com/chris-mount-60/sticky-bm-bridges}. The implementation is a short Python module. It uses floating point and computes the local time in
Algorithm~\ref{alg:bridge} by bisection. Figure~\ref{fig:lt-check} compares one component of Algorithm~\ref{alg:bridge} with its exact law. We take both ends
of the bridge at $0$, so that the whole bridge is the middle bridge of local time $r$ and operational length
$T-\stick r$, and the local time at a time $s$ has the explicit distribution function of
Proposition~\ref{prop:middle-lt}, which needs no numerical integration. With $r=0.8$, $T=2$, $\stick=0.5$ and
$s=0.9$, the figure shows this distribution function and the empirical distribution function of the local time
returned by $n=10^5$ independent runs of the algorithm, that is, the fraction of the returned values below each level.
On the whole range the two functions are indistinguishable. Their largest difference is $0.0013$, against the scale
$n^{-1/2}=0.0032$ of the fluctuations of an empirical distribution function, and the right panel magnifies the
neighbourhood where it occurs. Figure~\ref{fig:refined} shows a sample of the same implementation after $3000$ queries.

\FloatBarrier

\appendix

\section{Conditional laws and the Markov property}\label{app:A}

\subsection{Conditional laws}\label{app:conditional}

The samplers condition on values they have drawn. The events they condition on, such as $\{\Y_T=y,\ \locT_T=r\}$, have
probability zero, so we fix once which version of a conditional law we mean. Let $V$ and $\mathfrak w$ be the observed and conditioning vectors,
respectively. We specify their joint density in the free coordinates.
For observations at successive times, the local-time coordinates are
the gains between those times. Each gain has reference measure
$\delta_0(\md a)+\md a$, as in Proposition~\ref{prop:kernel}.
A recorded zero contributes its time as a free coordinate. Its
position is zero and its local time is that of a neighbouring
observation, as in Proposition~\ref{prop:zeros}.

When we condition on $\mathfrak w=(\Y_T,\locT_T)=(y,r)$, the gains
sum to $r$. The reference measure for the remaining coordinates can
therefore depend on $r$. For example, with one interior observation,
a zero gain after the observation gives the atom at $a=r$ in
\eqref{eq:bridge-one}. Proposition~\ref{prop:bridge-path}(a) gives
the reference measure for any finite set of observation times. A density is determined only almost
everywhere, so we name the versions. We always use the explicit formulas of Propositions~\ref{prop:kernel},
\ref{prop:zeros}, \ref{prop:middle} and~\ref{prop:record}, and their products under the Markov property. The
conditional law of $V$ given $\mathfrak w=w$ is then the ratio of densities
\begin{equation}\label{eq:bridge-def}
	\Probx{x}\{V\in\md v\mid \mathfrak w=w\}~\ass~\frac{\Probx{x}\{V\in\md v,\ \mathfrak w\in\md w\}}{\Probx{x}\{\mathfrak w\in\md w\}}
\end{equation}
at every $w$ at which the density of $\mathfrak w$ is positive. Three facts stand behind this convention. First, for any
choice of the two densities, \eqref{eq:bridge-def} is a regular conditional law of $V$ given $\mathfrak w$, for almost every
$w$, by Fubini's theorem. Second, for the explicit formulas we use, the ratio is a probability law at every $w$ at
which the denominator is positive, because the denominator is the integral of the numerator over $v$ at that $w$, by
the identities \eqref{eq:h-conv} and \eqref{eq:entrance} and the propositions cited, not only for almost every $w$.
Third, Proposition~\ref{prop:bridge-path} identifies the ratios for the endpoint $\mathfrak w=(\Y_T,\locT_T)$ with the
finite-dimensional laws of one measurable law of continuous paths, defined at every endpoint. The samplers only ever
condition on values at which the second fact applies.

The ratios \eqref{eq:bridge-def} are laws of finitely many values. The main case is the endpoint $\mathfrak w=(\Y_T,\locT_T)$,
which gives the bridges of Section~\ref{sec:bridge}. For it we also need the law of paths just mentioned, so that the
path can be evaluated at random times. Appendix~\ref{app:path-record} constructs this law. The position $\mathfrak w=\Y_T$
alone gives the bridges of $\Y$ treated in the companion note \citep{mountFatkullinDriftNote}, and $\mathfrak w$ may also contain zeros of $\Y$.

\subsection{The Markov property through the time change}\label{app:markov}

Section~\ref{sec:markov} uses the Markov property \eqref{eq:markov} at four kinds of times. These are fixed times, the
first zero $d_u$ after a fixed time $u$, the time $\hT_a+\stick a$ at which the local time first reaches a fixed level
$a>0$, and every finite $\filt_{\sigma_S}$-measurable time $S'\geq S$. The last kind includes a time chosen from the
values of $(\Y,\locT)$ drawn up to time $S$. All four are finite almost surely, since $\BM$ reaches every level, so that $\hT_a<\infty$, and $\BMr$ returns
to $0$ after every time, so that $d_u<\infty$. The next lemma transfers the
strong Markov property of $\BM$ through the time change.

\begin{lemma}\label{lem:markov}
Let $S\geq0$ be an almost surely finite random time such that $\sigma_S$ is a stopping time of $\BM$. Then
\eqref{eq:markov} holds. Each of the four kinds of times above has this property.
\end{lemma}

\begin{proof}
The clock $A_t=t+\stick\BMM_t$ is a continuous additive functional of $\BMr$, because $\BMM$ is the local time of
$\BMr$. A strong Markov process run on the inverse of such a functional is again strong Markov
\citep{volkonskiiRandomSubstitutionTime1958}; see \citet[Ch.~X, Thm.~(2.18)]{RevuzYor99} for a modern account for Brownian motion.
Concretely, $\BM'_t\ass\BM_{\sigma_S+t}-\BMM_{\sigma_S}$ is, given $\filt_{\sigma_S}$, a Brownian motion from $-\Y_S$,
and the process $(\Y',\locT')$ built from it as in Section~\ref{sec:time-change} satisfies $\Y_{S+v}=\Y'_v$ and
$\locT_{S+v}=\locT_S+\locT'_v$. This is \eqref{eq:markov}.

It remains to check that $\sigma_S$ is a stopping time for the times listed. Recall that $(\filt_t)$ is the
filtration generated by $\BM$, without completion. The strong Markov property of Brownian motion holds at its stopping times
(\citealp{huntTheoremsConcerningBrownian1956};\ \citealp[Thm.~2.6.15]{karatzas1991brownian}), and the first hitting time of a closed set by a continuous adapted process
after a stopping time is again one of them. For a fixed time $u$, the event
$\{\sigma_u\leq t\}=\{A_t\geq u\}$ lies in $\filt_t$. For the first zero after $u$, the time
$\sigma_{d_u}=\inf\{t\geq\sigma_u:\BMr_t=0\}$ is the first hitting time of $0$ by $\BMr$ after a stopping time. For the
local time, $\sigma_{\hT_a+\stick a}=\hT_a$. For an $\filt_{\sigma_S}$-measurable time $S'\geq S$, we have
$\{\sigma_{S'}\leq t\}=\{\sigma_S\leq t\}\cap\{A_t\geq S'\}$. This event lies in $\filt_t$, because $S'\id_{\{\sigma_S\leq t\}}$
is $\filt_t$-measurable by the definition of $\filt_{\sigma_S}$.
\end{proof}

\section{Brownian passage times and Bessel(3) bridges}\label{app:brownian}

This appendix proves the laws of Section~\ref{sec:passage} and gives a sampler for each. As there, $\BM$ starts at $0$
and $\hT_a$ is the first passage over the distance $a$.

\subsection{A passage time conditioned on a later one}

The law \eqref{eq:ptb} is the law at one time of the bridge of the first passage process. As a process in the level,
the first passage process is the right-continuous version $a\mapsto\hT_{a+}\ass\inf\{t:\BM_t>a\}$, the stable
subordinator of index $\frac12$, and $\hT_{a+}=\hT_a$ almost surely at each fixed level $a$. It is a first passage
bridge of Brownian motion read at one level \citep{bertoinPathTransformationsFirst2003}; the same bridges arise from
the inverse Gaussian subordinator, and \citet{ribeiroMonteCarloMethod2003} sample them for option pricing. The mixture
of two transformed inverse Gaussian variables below is the passage-time form of the sampler of
\citet[Prop.~2]{beskosRetrospectiveExactSimulation2006} for the time of the minimum of a Brownian bridge. We include
its derivation in our notation. We write $\IG(m,\lambda)$ for the inverse Gaussian law with mean $m>0$
and shape $\lambda>0$. Its density is
\[
	\sqrt{\frac{\lambda}{2\mpi v^3}}\,\exp\Big\{-\frac{\lambda(v-m)^2}{2m^2v}\Big\},\qquad v>0,
\]
and one draw from it uses one standard normal and one uniform \citep{michaelGeneratingRandomVariates1976}.

\begin{proposition}\label{prop:ptb}
Let $a,b,t>0$, let $\uni$ be uniform on $(0,1)$, and let $V\law\IG(a/b,\,a^2/t)$ and $X\law\IG(b/a,\,b^2/t)$, all
independent. Then
\[
	\Probx{0}\big\{\hT_a\in\md u\mid\hT_{a+b}=t\big\}~=~
	\Prob\Big\{\id\Big\{\uni<\tfrac{b}{a+b}\Big\}\frac{tV}{1+V}+\id\Big\{\uni\geq\tfrac{b}{a+b}\Big\}\frac{t}{1+X}\in\md u\Big\}.
\]
\end{proposition}

\begin{proof}
Write $\psi(u)$ for the density in \eqref{eq:ptb}, and change variables to $v=u/(t-u)$. Then $u=tv/(1+v)$,
$t-u=t/(1+v)$ and $\md u=t(1+v)^{-2}\md v$. Since $a^2+b^2-(a+b)^2=-2ab$ and $a^2/v+b^2v-2ab=(bv-a)^2/v$,
\[
	\psi\Big(\frac{tv}{1+v}\Big)\frac{t}{(1+v)^2}~=~\frac{ab}{(a+b)\sqrt{2\mpi t}}\,
	\big(v^{-3/2}+v^{-1/2}\big)\exp\Big\{-\frac{(bv-a)^2}{2tv}\Big\},\qquad v>0 .
\]
The term with $v^{-3/2}$ is $b/(a+b)$ times the density of $\IG(a/b,a^2/t)$ at $v$. The term with $v^{-1/2}$ is
$a/(a+b)$ times the density at $v$ of the reciprocal of an $\IG(b/a,b^2/t)$ random variable, and
$u=tv/(1+v)=t/(1+1/v)$.
\end{proof}

The same change of variables $v=u/(t-u)$ gives the distribution function \eqref{eq:ptb-cdf}.

\subsection{Bessel(3) bridges}

We sample the value at a time $u$ of a path from $z$ to $z'$ over $[0,t]$ that does not reach $0$ in between, in the
three cases of Section~\ref{sec:passage}. Let $z,z'\geq0$ with $z+z'>0$, let $t>0$ and $0<u<t$, and let $\uni$, $\gauss$
and $\erv$ be independent, with $\uni$ uniform on $(0,1)$, $\gauss$ standard normal and $\erv$ exponential with mean
$1$. Put $\kappa\ass zz'/t$,
\[
	\eta~\ass~
	\begin{cases}
		1+\dfrac1\kappa\log\big(1-\uni\,(1-\me^{-2\kappa})\big), & \kappa>0,\\[4pt]
		1, & \kappa=0,
	\end{cases}
\]
and
\begin{equation}\label{eq:bessel-sampler}
	m^2~\ass~\Big(1-\frac ut\Big)^2z^2+\frac{u^2}{t^2}\,z'^2+2\,\frac ut\Big(1-\frac ut\Big)zz'\,\eta,\qquad
	\bessb_u~\ass~\Big(\Big(m+\sqrt{\tfrac{u(t-u)}t}\,\gauss\Big)^2+\frac{2u(t-u)}t\,\erv\Big)^{1/2},
\end{equation}
with $m\geq0$.

\begin{lemma}\label{lem:bessel}
Let $\bessb_u$ be given by \eqref{eq:bessel-sampler}. Then, for $w>0$,
\begin{align*}
	\Prob\{\bessb_u\in\md w\}~&=~\frac{p_{\mathfrak A}(w,u\mid z)\,p_{\mathfrak A}(z',t-u\mid w)}{p_{\mathfrak A}(z',t\mid z)}\md w
	&&z,z'>0,\\
	\Prob\{\bessb_u\in\md w\}~&=~\frac{h(w,u)\,p_{\mathfrak A}(z',t-u\mid w)}{h(z',t)}\md w&&z=0,\\
	\Prob\{\bessb_u\in\md w\}~&=~\frac{p_{\mathfrak A}(w,u\mid z)\,h(w,t-u)}{h(z,t)}\md w&&z'=0 .
\end{align*}
\end{lemma}

The first law is the law at time $u$ of Brownian motion from $z$ conditioned not to reach $0$ before time $t$
and to be at $z'$ at time $t$. The second and third are the laws after the last zero of a path and before its first
zero, and each is the time reversal of the other. When $zz'=0$, the uniform is not used, since $\eta$ enters $m$ only
through the product $zz'$. A Bessel(3) bridge is the length of a three-dimensional Brownian motion conditioned on the length of its endpoint.
Given the direction of the endpoint, it is the length of a three-dimensional Brownian bridge between two fixed points
\citep[Ch.~XI, Exercise~(3.8)]{RevuzYor99}, and the uniform samples the cosine $\eta$ of the angle between the two
endpoints, which is all that the law of the length needs. When one end is at $0$, the law does not depend on the
direction, and the exact algorithms for diffusions sample the bridge in this way
\citep{beskosRetrospectiveExactSimulation2006}.

\begin{proof}
Let $\boldsymbol{\mathfrak W}$ be a three-dimensional Brownian motion from $\mathbf z\ass(z,0,0)$. Its length $|\boldsymbol{\mathfrak W}|$ is a Bessel(3) process,
and the three displays are the laws at time $u$ of its bridges \citep[Ch.~XI, Def.~(3.1)]{RevuzYor99}. It remains to see that $\bessb_u$ has the law of $|\boldsymbol{\mathfrak W}_u|$ given $|\boldsymbol{\mathfrak W}_t|=z'$. Let first $z,z'>0$. On the sphere $\{|b|=z'\}$, the density of
$\boldsymbol{\mathfrak W}_t$ is proportional to $\me^{z b_1/t}=\me^{\kappa\eta'}$, where $\eta'=b_1/z'$ is the cosine of the angle between $b$ and
$\mathbf z$. Under the uniform law on the sphere $\eta'$ is uniform on $[-1,1]$, so given $|\boldsymbol{\mathfrak W}_t|=z'$ it has the density
$\kappa\me^{\kappa\eta'}/(2\sinh\kappa)$ on $[-1,1]$. Its distribution function at $\eta$ is $1-\uni$, so $\eta$ has this law.
Given $\boldsymbol{\mathfrak W}_t=b$, the process $\boldsymbol{\mathfrak W}$ is a three-dimensional Brownian bridge, and
\[
	\boldsymbol{\mathfrak W}_u~\law~\mathbf m+\sqrt{\frac{u(t-u)}t}\,(\gauss_1,\gauss_2,\gauss_3),\qquad \mathbf m~\ass~\Big(1-\frac ut\Big)\mathbf z+\frac ut\,b,
\]
with independent standard normals $\gauss_1,\gauss_2,\gauss_3$. The normal vector is invariant under rotations, so the
law of $|\boldsymbol{\mathfrak W}_u|$ depends on $\mathbf m$ only through $|\mathbf m|=m$. Rotating $\mathbf m$ to $(m,0,0)$ gives
\[
	|\boldsymbol{\mathfrak W}_u|^2~\law~\Big(m+\sqrt{\tfrac{u(t-u)}t}\,\gauss_1\Big)^2+\frac{u(t-u)}t\big(\gauss_2^2+\gauss_3^2\big),
\]
and $\gauss_2^2+\gauss_3^2$ has the law of $2\erv$. For $z=0$ or $z'=0$ the direction of the endpoint is immaterial,
and the same computation applies.
\end{proof}

\section{Gamma envelopes}\label{app:envelope}

The loop of Section~\ref{sec:middle} samples densities proportional to
$w^j\me^{-kw-w^2/2}$ on $w>0$, with $k\geq0$. For fixed $k$ and $\lambda>0$, the function $w^j\me^{-kw-w^2/2}$ is at most
$\me^{(\lambda-k)^2/2}w^j\me^{-\lambda w}$. The ratio of their masses is largest when the derivative of its logarithm in
$\lambda$, which is $(j+1)/\lambda-(\lambda-k)$, vanishes. This gives the rate $\lambda_j$ and the ratio $A_j(k)$ of
Section~\ref{sec:middle}.

\begin{lemma}\label{lem:envelope}
Let $k\geq0$ and let $j\geq0$ be an integer. For $w>0$,
\[
	w^j\me^{-kw-w^2/2}~=~\me^{(\lambda_j-k)^2/2}\,w^j\me^{-\lambda_jw}\cdot\exp\Big\{-\frac{(w-\lambda_j+k)^2}2\Big\}.
\]
The ratio $A_j(k)$ is strictly increasing in $k\in[0,\infty)$ and tends to $1$ as $k\to\infty$, and for $j\leq2$,
\[
	A_j(k)~\geq~A_2(0)~=~\frac{3\sqrt3}2\sqrt{\frac\mpi2}\,\me^{-3/2}~=~0.72656\ldots
\]
\end{lemma}

By the identity, $w^j\me^{-kw-w^2/2}$ is at most $\me^{(\lambda_j-k)^2/2}j!\,\lambda_j^{-j-1}$ times the density of
$\GammaLaw(j+1,\lambda_j)$.

\begin{proof}
The identity is
\[
	-kw-\frac{w^2}2~=~\frac{(\lambda-k)^2}2-\lambda w-\frac{(w-\lambda+k)^2}2
\]
at $\lambda=\lambda_j$. With
$I_m\ass\int_0^\infty w^m\me^{-kw-w^2/2}\md w$ as in Section~\ref{sec:middle}, we have $\frac{\md}{\md k}I_m=-I_{m+1}$. The rate $\lambda_j$ is the
positive root of $\lambda(\lambda-k)=j+1$, so $\lambda_j-k=(j+1)/\lambda_j$. In the derivative of $\log A_j$, the
terms that come from $\lambda_j'(k)$ cancel,
\[
	\frac{\md}{\md k}\log A_j(k)~=~\Big(\frac{j+1}{\lambda_j}-(\lambda_j-k)\Big)\lambda_j'(k)+(\lambda_j-k)-\frac{I_{j+1}}{I_j}
	~=~(\lambda_j-k)-\frac{I_{j+1}}{I_j}.
\]
Hence $A_j$ increases exactly when the mean $I_{j+1}/I_j$ of the density proportional to $w^j\me^{-kw-w^2/2}$ is below
$\lambda_j-k$, the mean of the Gamma proposal. The derivative of $w^{j+1}\me^{-kw-w^2/2}$ integrates to $0$, because the
function vanishes at $w=0$ and as $w\to\infty$. This gives $0=(j+1)I_j-kI_{j+1}-I_{j+2}$. By the Cauchy--Schwarz
inequality, $I_{j+1}^2<I_jI_{j+2}$, so with $m\ass I_{j+1}/I_j$,
\[
	j+1~=~k\,m+\frac{I_{j+2}}{I_j}~>~k\,m+m^2 .
\]
Thus $m$ is a positive number at which $q\mapsto q^2+kq-(j+1)$ is negative. Hence $m$ is below the positive root of this
quadratic, which is $\lambda_j-k$, and $A_j$ is strictly increasing. As $k\to\infty$, the substitution $w=w'/k$ gives
\[
	A_j(k)~=~\Big(\frac{\lambda_j}k\Big)^{j+1}\me^{-(\lambda_j-k)^2/2}\,\frac1{j!}\int_0^\infty w'^j\me^{-w'-w'^2/2k^2}\md w'
	~\longrightarrow~1,
\]
since $\lambda_j-k=(j+1)/\lambda_j\to0$. At $k=0$ we have $\lambda_j=\sqrt{j+1}$, $I_0=I_2=\sqrt{\mpi/2}$ and $I_1=1$.
This gives $A_0(0)=\sqrt{\mpi/2}\,\me^{-1/2}=0.760$, $A_1(0)=2\me^{-1}=0.736$ and $A_2(0)=0.727$.
\end{proof}

\section{The path given a skeleton}\label{app:path-record}

The results of Section~\ref{sec:record-law} need the conditional law of the whole path given a skeleton, as one kernel
that is defined at every skeleton, and they need the skeletons of a single path to behave like the queries. This
appendix sets up both, and it finds the law of a skeleton under the kernel. The kernel uses the Bessel(3) bridges and the
middle bridge as laws of paths, and Appendix~\ref{app:bridges-paths} defines them.

We work on the space $\mathcal C$ of continuous functions from $[0,\infty)$ to $[0,\infty)^2$, with the topology of
uniform convergence on compact sets, and on the space $\mathcal C[0,T]$ of continuous functions on $[0,T]$, with the
uniform topology. We write $(\Y,\locT)$ also for the coordinate process, so that $\Y_t(w)$ and $\locT_t(w)$ are the
position and the local time of the path $w$ at time $t$. The Borel $\sigma$-field of $\mathcal C$ is generated by the
coordinates \citep[Ch.~XIII, Prop.~(1.2)]{RevuzYor99}. Two laws on $\mathcal C$ are therefore equal when their
finite-dimensional laws are equal. The time change of Section~\ref{sec:time-change} is a measurable map from
$\Omega$ to $\mathcal C$ that takes the Brownian path to $(\Y,\locT)$, and we write $\Probx{x}$ also for the image
of $\Probx{x}$ under it, the law of $(\Y,\locT)$ on $\mathcal C$. For a path $w$
and a finite set $\mathcal U\ni0$ of times, $\rrec(\mathcal U)(w)$ is the set \eqref{eq:revealed} computed from $w$, with the zeros $g_u$
and $d_u$ computed from $\Y(w)$.

\subsection{The skeletons of a path}\label{app:revealed}

A query adds a point to a skeleton, and the theorem compares this with adding a
time to $\mathcal T$. Both operations must give the same set of points at every time, not only at almost every one. This is a
statement about a single path, and it needs three of its properties. We write $\Paths$ for the set of paths
$w\in\mathcal C$ with the following properties.
\begin{enumerate}
    \renewcommand{\theenumi}{P\arabic{enumi}}
    \renewcommand{\labelenumi}{(\theenumi)}
    \item\label{eq:P1} The local time $\locT$ starts at $0$, is nondecreasing, and is constant on every interval on which $\Y>0$.
    \item\label{eq:P2} The zero set of $\Y$ has no isolated points in $(0,\infty)$.
    \item\label{eq:P3} Every interval on which $\locT$ is constant contains at most two zeros of $\Y$.
\end{enumerate}
On $\mathcal C[0,T]$ we use the same properties with $(0,T)$ in place of $(0,\infty)$. For a path $w$ that
satisfies \eqref{eq:P1}, the set $\rrec(\mathcal U)(w)$ is a skeleton, since the local time is constant on $[u,d_u]$ and on $[g_{u'},u']$, and
$d_u<g_{u'}$ when $\locT_u<\locT_{u'}$.

\begin{lemma}\label{lem:path-records}
Let $w\in\Paths$, let $\mathcal U\ni0$ be a finite set of times, and put $\rec\ass\rrec(\mathcal U)(w)$. Let $u$ be a time.
\begin{enumerate}
	\item[(a)] If $u$ is the time of a point of $\rec$, then $\rrec(\mathcal U\cup\{u\})(w)=\rec$.
	\item[(b)] Otherwise $\rrec(\mathcal U\cup\{u\})(w)=\rec\cup\{(u,\Y_u,\locT_u)\}\cup\mathcal Z$, where the set $\mathcal Z$
	of new recorded zeros is, with $(s_i,z_i,a_i)$, $0\leq i\leq m$, the points of $\rec$,
	\[
		\mathcal Z~=~
		\begin{cases}
			\{(g_u,0,\locT_u),\,(d_u,0,\locT_u)\}, & s_i<u<s_{i+1},\ \delta_i>0,\ \Y_u>0,\\[2pt]
			\{(d_{s_m},0,a_m):z_m>0\}\cup\{(g_u,0,\locT_u):\Y_u>0\}, & u>s_m,\ \locT_u>a_m,\\[2pt]
			\emptyset & \text{otherwise.}
		\end{cases}
	\]
	\item[(c)] If $\mathcal U\subset \mathcal U'$, then $\rrec(\mathcal U)(w)\subset\rrec(\mathcal U')(w)$. If $\mathcal Z$ is a finite set of times of points of
	$\rrec(\mathcal U')(w)$, then $\rrec(\mathcal U'\cup \mathcal Z)(w)=\rrec(\mathcal U')(w)$.
\end{enumerate}
\end{lemma}

The cases in (b) are those of the query in Section~\ref{sec:records}. The skeleton of the path at $\mathcal U\cup\{u\}$ is the
skeleton at $\mathcal U$ with the points that the query at $u$ adds, computed from the path.

\begin{proof}
Only the consecutive times $s<s'$ of $\mathcal U$ around $u$ matter, or the last time $s$ of $\mathcal U$ if $u$ is after it. The skeleton
$\rrec(\mathcal U\cup\{u\})(w)$ differs from $\rec$ only through the pairs $(s,u)$ and $(u,s')$.

(a) If $u\in \mathcal U$ there is nothing to prove. Let $u=d_s$ be a recorded zero, so that $\Y_s>0$ and $\locT_s<\locT_{s'}$.
By \eqref{eq:P1} the local time is constant on $[s,d_s]$, so the pair $(s,d_s)$ adds no zero. The pair $(d_s,s')$ adds
$g_{s'}$ if $\Y_{s'}>0$, and no first zero, since $\Y_{d_s}=0$. This is $\rec$. A recorded zero $g_{s'}$ is the same
case with time reversed.

(b) Let $\locT_s=\locT_{s'}$. Then $\locT_u=\locT_s$, both pairs are constant, and $\rrec(\mathcal U\cup\{u\})(w)$ adds only the
point at $u$. Now let $\locT_s<\locT_{s'}$, and let $d\ass d_s$ if $\Y_s>0$ and $d\ass s$ otherwise, and
$g\ass g_{s'}$ if $\Y_{s'}>0$ and $g\ass s'$ otherwise. If $u\in(s,d)$, then $\Y>0$ on $[s,d)$, the pair $(s,u)$ is
constant, and the pair $(u,s')$ adds $d_u=d$ and $g_{s'}$, which are already in $\rec$. The case $u\in(g,s')$ is the
same with time reversed. It remains to take $u$ in the increasing interval $(d,g)$, and we treat the pair $(s,u)$.
The pair $(u,s')$ is the same case with time reversed.

If $\locT_s<\locT_u$, the pair $(s,u)$ adds $d_s$, which is $d\in\rec$ when $\Y_s>0$, and $g_u$ when $\Y_u>0$. Let
$\locT_s=\locT_u$, so that $\locT$ is constant on $[s,u]$. If $\Y_s>0$, then $d$ is a zero with $\Y>0$ on $(s,d)$. By
\eqref{eq:P2} zeros of $\Y$ accumulate at $d$ from the right, so $[s,u]$ contains infinitely many zeros, against \eqref{eq:P3}. Hence
$\Y_s=0$ and $d=s$. If $\Y_u>0$ and $g_u>s$, then by \eqref{eq:P2} zeros accumulate at $g_u$ from the left, and \eqref{eq:P3} fails in the
same way. Hence $g_u=s$, and the point $(g_u,0,\locT_u)=(s,0,\locT_s)$ is already in $\rec$. In every case the pair adds
the points listed in (b). After the last time $s$ of $\mathcal U$, the pair $(s,u)$ adds exactly the points listed in (b), by the
definition \eqref{eq:revealed}.

(c) By (b), adding one time to $\mathcal U$ only adds points, and the first claim follows by adding the times of $\mathcal U'\setminus \mathcal U$
one at a time. The second claim follows from (a) in the same way.
\end{proof}

\begin{lemma}\label{lem:paths-full}
For every $x\geq0$, $\Probx{x}(\Paths)=1$.
\end{lemma}

\begin{proof}
Property \eqref{eq:P1} holds by the construction of Section~\ref{sec:time-change}. The inverse clock $\sigma$ is continuous and
strictly increasing from $[0,\infty)$ onto $[0,\infty)$. It maps the zeros of $\Y$ onto the zeros of $\BMr$, and the
intervals on which $\locT$ is constant onto the intervals on which $\BMM$ is constant. It is therefore enough to prove
\eqref{eq:P2} and \eqref{eq:P3} for $\BMr$ and $\BMM$. Before $\hT_0$ we have $\BMr=-\BM>0$. By the strong Markov property at $\hT_0$
\citep{huntTheoremsConcerningBrownian1956}, $\BM_{\hT_0+\iv}$ is a Brownian motion from $0$, and Lévy's theorem
\citep{levyProcessusStochastiquesMouvement1965} makes $\BMr_{\hT_0+\iv}$ a reflected Brownian motion. The zero set of a Brownian motion has no isolated
points, and $0$ is a limit of its zeros \citep[Thm.~2.9.6]{karatzas1991brownian}. This is \eqref{eq:P2}.

For \eqref{eq:P3}, let $\BMM$ be constant, equal to $m$, on an interval that contains three zeros $t_1<t_2<t_3$ of $\BMr$. Then
$\BM_{t_1}=\BM_{t_2}=m$, and $\BM\leq m$ on $[0,t_3]$. So $\BM$ has local maxima at $t_1$ and at $t_2$ with the same
value $m$. Almost surely no two local maxima of $\BM$ have the same value. Indeed, each local maximum is the maximum
over an interval with rational ends, and for rationals $q_1<q_2<q_3<q_4$ the maximum over $[q_3,q_4]$ is $\BM_{q_3}$
plus a variable that is independent of $(\BM_v)_{v\leq q_3}$ and has a density, so it differs almost surely from the
maximum over $[q_1,q_2]$.
\end{proof}

\subsection{Bessel(3) and middle bridges as laws of paths}\label{app:bridges-paths}

The kernel of Appendix~\ref{app:kernel} needs the Bessel(3) bridges and the middle bridge as laws of paths, not only
at one time. The facts are classical. The Bessel(3) process is the Doob $\mathfrak h$-transform of Brownian motion
killed at $0$ \citep{doobConditionalBrownianMotion1957a}, and Brownian motion killed at $0$ and conditioned on its
endpoint, or on its first passage time to $0$, is a Bessel(3) bridge
\citep{williamsPathDecompositionContinuity1974,pitmanOnedimensionalBrownianMotion1975}. See \citet[Ch.~XI, \S3]{RevuzYor99} for a modern account. A Brownian motion
conditioned on a first passage time is a first passage bridge \citep{bertoinPathTransformationsFirst2003}. We state
these facts in the form that we use. The laws below live on spaces $\mathcal C[0,t]$ whose length $t$ varies. For the
statements about kernels, we regard each of them as a law on $\mathcal C$ by extending every path beyond its interval
by its final value; a law on $\mathcal C[0,t]$ and its extension determine each other.

Let $z,z'\geq0$ with $z+z'>0$, and $t>0$. Let $\boldsymbol{\mathfrak W}$ be a three-dimensional Brownian motion from $(z,0,0)$, and
condition it on $|\boldsymbol{\mathfrak W}_t|=z'$ as in the proof of Lemma~\ref{lem:bessel}. That is, the endpoint $\boldsymbol{\mathfrak W}_t$ is a
point of the sphere of radius $z'$ with the direction law of that proof, and given $\boldsymbol{\mathfrak W}_t=b$ the process
$\boldsymbol{\mathfrak W}$ is the three-dimensional Brownian bridge from $(z,0,0)$ to $b$ over $[0,t]$. We set
\[
    \bessb_v\ass|\boldsymbol{\mathfrak W}_v|,\qquad 0\leq v\leq t.
\]
The pair $(\bessb_v,0)_{0\leq v\leq t}$ has law $\Bes_{z,z',t}$ on $\mathcal C[0,t]$.

To construct the middle bridge, let $\delta,t>0$, take $(\bessb,0)$ with law $\Bes_{\delta,0,t}$, and
put $\BMb_v\ass\delta-\bessb_v$ for $v\leq t$. Then $\BMb$ starts at $0$ and reaches the level $\delta$ for the first
time at $t$. We build its running maximum $\BMM$, the reflected process $\BMr=\BMM-\BMb$, the clock $A$ and its inverse
$\sigma$ from $\BMb$ on $[0,t]$ as in Section~\ref{sec:time-change}. Then $A_t=t+\stick\delta$.
We take $\Mid_{\delta,t}$ to be the law on $\mathcal C[0,t+\stick\delta]$ of the path
\[
    (\BMr_{\sigma_v},\BMM_{\sigma_v})_{0\leq v\leq t+\stick\delta}.
\]

\begin{lemma}\label{lem:bridges-paths}
\begin{enumerate}
	\item[(a)] Let $0=v_0<v_1<\dots<v_n<v_{n+1}=t$, and put $z_0\ass z$ and $z_{n+1}\ass z'$. Under $\Bes_{z,z',t}$,
	\[
		\Bes_{z,z',t}\big\{\bessb_{v_1}\in\md z_1,\dots,\bessb_{v_n}\in\md z_n\big\}
		~=~\frac{\prod_{j=0}^{n}k_j}{k}\,\md z_1\cdots\md z_n,\qquad z_1,\dots,z_n>0,
	\]
	where $k_j\ass p_{\mathfrak A}(z_{j+1},v_{j+1}-v_j\mid z_j)$, except that $k_0\ass2h(z_1,v_1)$ if $z=0$ and
	$k_n\ass h(z_n,t-v_n)$ if $z'=0$, and $k$ is the same expression for $n=0$. The path $\bessb$ is positive on $(0,t)$.
	The law $\Bes_{0,z',t}$ is the law of the time reversal $(\bessb_{t-v})_{v\leq t}$ under $\Bes_{z',0,t}$, and
	$(z,z',t)\mapsto\Bes_{z,z',t}$ is a kernel.
	\item[(b)] Let $z,t>0$. Under $\Probx{z}$, the process $-\BM$ is a Brownian motion from $z$, and $\hT_0$ is its first
	passage time to $0$. Let $\varsigma$ be a stopping time of the coordinate filtration of paths on $[0,t]$, with
	$\Bes_{z,0,t}\{\varsigma<t\}=1$. For every nonnegative measurable functional $F$ of the path stopped at $\varsigma$,
	\[
		\Bes_{z,0,t}\big[F(\bessb_{\iv\wedge\varsigma})\big]
		~=~\frac{\E_z\big[F(-\BM_{\iv\wedge\varsigma})\,\id_{\{\varsigma<\hT_0\}}\,h(-\BM_\varsigma,\ t-\varsigma)\big]}{h(z,t)} ,
	\]
	with $\varsigma$ evaluated on $\bessb$ on the left and on $-\BM$ on the right.
	\item[(c)] Let $\delta,t>0$ and $\Delta\ass t+\stick\delta$, and let $\varsigma$ be a stopping time as in (b), with
	$\varsigma<t$ for the path $\BMb$. For every nonnegative measurable functional $F$ of the path stopped at
	$\varsigma$,
	\[
		\E\big[F(\BMb_{\iv\wedge\varsigma})\big]
		~=~\frac{\E_0\big[F(\BM_{\iv\wedge\varsigma})\,\id_{\{\varsigma<\hT_\delta\}}\,h(\delta-\BM_\varsigma,\ t-\varsigma)\big]}{h(\delta,t)} ,
	\]
	where $\E$ is the expectation of the construction of $\BMb$. In particular $\Mid_{\delta,t}$ is the middle bridge of
	Section~\ref{sec:bridge-given-zeros}, and on the physical times before $A_{t'}$, for every $t'<t$, it is absolutely continuous
	with respect to $\Probx{0}$. Under $\Mid_{\delta,t}$ the local time satisfies $\locT_v<\delta=\locT_\Delta$ for
	$v<\Delta$, and it increases on every neighbourhood of $0$ and of $\Delta$ inside $[0,\Delta]$. The map
	$(\delta,t)\mapsto\Mid_{\delta,t}$ is a kernel.
\end{enumerate}
\end{lemma}

\begin{proof}
(a) The display gives the finite-dimensional laws of the Bessel(3) bridge \citep[Ch.~XI, Def.~(3.1)]{RevuzYor99},
and reversal follows from $p_{\mathfrak A}(a,s\mid b)=p_{\mathfrak A}(b,s\mid a)$. The path is continuous since
$\boldsymbol{\mathfrak W}$ is, and positive on $(0,t)$ since three-dimensional Brownian motion does not hit the origin after time $0$.
The sets $B$ for which $\Bes_{z,z',t}(B)$ is measurable in $(z,z',t)$ form a $\lambda$-system that contains the
cylinder sets, so $\Bes$ is a kernel by the $\pi$-$\lambda$ theorem.

(b) This is Doob's $\mathfrak h$-transform \citep{doobConditionalBrownianMotion1957a} at a stopping time. By (a), the
display holds at every fixed time $\varsigma=s<t$. The process $M_s\ass\id_{\{s<\hT_0\}}\,h(-\BM_s,t-s)$ is a bounded
martingale under $\Probx{z}$ on $[0,t']$ for every $t'<t$, by the Markov property, the symmetry of $p_{\mathfrak A}$
and \eqref{eq:entrance}. Fix $t'<t$. The functional $F\cdot\id_{\{\varsigma\leq t'\}}$ is a functional of the path stopped
at the bounded stopping time $\varsigma\wedge t'$, and optional sampling at $\varsigma\wedge t'$ gives the display with
$\id_{\{\varsigma\leq t'\}}$ inside both expectations. Both sides increase as $t'\uparrow t$. Monotone convergence gives
(b), since $\Bes_{z,0,t}\{\varsigma<t\}=1$ on the left and $h(\cdot,t-\varsigma)=0$ for $\varsigma\geq t$ on the right.

(c) Under $\Probx{\delta}$ the process $-\BM$ is a Brownian motion from $\delta$, and $\delta-(-\BM)$ is a Brownian
motion from $0$ whose first passage to $\delta$ is $\hT_0$. So (c) is (b) with $z=\delta$, read for $\BMb=\delta-\bessb$. For a fixed time $\varsigma=s$, the ratio $h(\delta-\BM_s,t-s)/h(\delta,t)$ on $\{s<\hT_\delta\}$ is
the conditional density of $\hT_\delta$ at $t$ given $\filt_s$, divided by its density. So the law of $\BMb$ is the
conditional law of $\BM$ given $\hT_\delta=t$, as a ratio of densities, and $\Mid_{\delta,t}$ is the middle bridge. On
$[0,t']$ the density of the law of $\BMb$ with respect to the law of $\BM$ is the bounded variable
$M_{t'}/h(\delta,t)$, and the path before the physical time $A_{t'}$ is a function of $\BMb$ on $[0,t']$. The running
maximum of $\BMb$ is below $\delta$ before $t$ and equals $\delta$ at $t$. It therefore increases on every neighbourhood
of $t$, and the local time increases on every neighbourhood of $\Delta$. Near $0$ the same holds under $\Probx{0}$, and
absolute continuity transfers it. Finally, $\Mid_{\delta,t}$ is the image of the kernel $\Bes_{\delta,0,t}$ under a
map that is measurable in $(\delta,t,\bessb)$, so it is a kernel.
\end{proof}

\subsection{The conditional path kernel}\label{app:kernel}

We now fix the skeletons on which the kernel is defined. We write $\Recs$ for the set of
skeletons $\rec$ with the following properties. The first point is $(0,z_0,0)$ for some $z_0\geq0$, the times of the points are distinct, and
the local times are nondecreasing in time. Each constant interval has $z_i+z_{i+1}>0$. Each increasing interval has
$\Delta_i>\stick\delta_i$. Every point at position $0$ is an end of an increasing interval, except possibly the last
point. We write $\Recs_T$ for the skeletons in $\Recs$ whose last point is at the time $T$, and whose last point, if at
position $0$, ends an increasing interval. The set $\Recs$ is a Borel subset of the union over $m$ of the spaces
$([0,\infty)^3)^{m+1}$ of ordered points. For $\rec\in\Recs$ we write $\mathcal T_\rec$ for the set of times of its points.

Let $\rec\in\Recs$ have points $(s_i,z_i,a_i)$, $0\leq i\leq m$. To construct $\Ker_\rec$, take independent bridges
with the interval laws of Section~\ref{sec:record-law}. On interval $i$ shift the time origin to $s_i$ and add $a_i$
to the local-time coordinate. After $s_m$ append an independent path of law $\Probx{z_m}$ with initial local time
shifted by $a_m$. For $\rec\in\Recs_T$ the construction ends at $T$ and gives a law on $\mathcal C[0,T]$.
The shifted paths start and end at the values that $\rec$ prescribes, so the path is continuous almost surely. On the
null set where two of them do not agree at a common time, we set the path to $(0,0)$.

For a path $w$ of $\Ker_\rec$ we have $\rrec(\mathcal T_\rec)(w)=\rec$ almost surely. The local time of $w$ is constant on the
constant intervals of $\rec$ and increases on the increasing ones, whose ends are at $0$, so no pair of consecutive
times of $\mathcal T_\rec$ records a zero. By Lemma~\ref{lem:path-records}(b), $\rec_u(w)$ is then the skeleton that the query at
$u$ would return if it read its new points from $w$. The query is a measurable function of $u$, of $\rec$ and of the elementary random variables it draws, so
$\qry_u(\rec,\iv)$ is a kernel in $(\rec,u)$.

\begin{lemma}\label{lem:kernel}
The map $\rec\mapsto\Ker_\rec$ is a probability kernel from $\Recs$ to $\mathcal C$, and from $\Recs_T$ to
$\mathcal C[0,T]$. For every $\rec\in\Recs$ or $\Recs_T$, $\Ker_\rec(\Paths)=1$. The map $(\rec,u,w)\mapsto\rec_u(w)$
is measurable.
\end{lemma}

\begin{proof}
By Lemma~\ref{lem:bridges-paths}, $\Bes_{z,z',\Delta}$ and $\Mid_{\delta,t}$ are kernels in their parameters, and so is
$\Probx{z}$, the image of the law of $\BM$ under a map that is measurable in $(z,\BM)$. The path is a measurable
function of $\rec$ and of the shifted paths. The image of a product of kernels under a measurable map is again a kernel
\citep[Lemmas~3.2 and~3.3]{kallenbergFoundationsModernProbability2021}.

Property \eqref{eq:P1} holds on each shifted path, and so on the whole path. For \eqref{eq:P2}, the zeros of the path lie inside the
middle bridges, after the last point, or at points of $\rec$ at position $0$. Inside a middle bridge and after the last
point there are no isolated zeros, by Lemma~\ref{lem:paths-full} and the absolute continuity of
Lemma~\ref{lem:bridges-paths}(c). In a middle bridge the local time increases on every neighbourhood of each end
inside it (Lemma~\ref{lem:bridges-paths}(c)). So does the local time from $0$ under $\Probx{0}$, since
$\Probx{0}\{\locT_s=0\}=0$ for $s>0$ by Proposition~\ref{prop:kernel}. By the definition of $\Recs$, every point at
position $0$ is therefore a limit of zeros. For \eqref{eq:P3}, let $I$ be an interval on which the local time is constant. By
the same fact, $I$ meets a middle bridge, or the path from a last point at $0$, only inside one constant interval of
that path, where \eqref{eq:P3} holds by Lemmas~\ref{lem:paths-full} and~\ref{lem:bridges-paths}. Otherwise $I$ runs through
consecutive Bessel(3) bridges, which have no zeros inside. The points that join two of them are positive, by the
definition of $\Recs$. So the zeros in $I$ are at most the two ends of this run, or one end and the first zero after a
positive last point, after which the local time increases at once.

The skeleton $\rec_u(w)$ is built from $\rec$, from $(\Y_u(w),\locT_u(w))$ and from $g_u(w)$ and $d_u(w)$. The last two
are measurable in $(u,w)$. Indeed, for $0<c\leq u$ the event $\{g_u<c\}$ is the event that $\Y$ is positive on
$[c,u]$, and for $c\geq u$ the event $\{d_u>c\}$ is the event that $\Y$ is positive on $[u,c]$; both are measurable
in $(u,c,w)$, since $\Y$ is continuous and it suffices to test rational times.
\end{proof}

\subsection{The law of a skeleton}\label{app:record-law}

The proofs of Appendix~\ref{app:refine-proof} use the law of the skeleton $\rec_{\mathcal V}$ under $\Ker_\rec$, which is a
product over the intervals. A time at which the position is $0$ plays one of two roles in a density. It is either a
prescribed time, at which the law of the state has an atom at position $0$, or a recorded zero, whose time has a
density. The two roles give different factors for the intervals that end there, and the factor of an interval depends
on the role of each of its ends. We say that an end is read as a zero in the second role and as fixed in the first.
Write $\mathcal O$ for the set of ends of intervals of $\rec_{\mathcal V}$ that are read as zeros. For the interval
$(s_i,s_{i+1})$ of $\rec_{\mathcal V}$, the roles of its ends determine the factor $\varphi_i$ as follows.
\[
	\varphi_i~\ass~
	\begin{cases}
		p_{\mathfrak A}(z_{i+1},\Delta_i\mid z_i), & \delta_i=0,\ s_i,s_{i+1}\notin\mathcal O,\\
		h(z_i,\Delta_i), & \delta_i=0,\ s_{i+1}\in\mathcal O,\\
		2h(z_{i+1},\Delta_i), & \delta_i=0,\ s_i\in\mathcal O,\\
		h(\delta_i,\Delta_i-\stick\delta_i), & \delta_i>0,\ s_{i+1}\notin\mathcal V,\\
		\stick\,h(\delta_i,\Delta_i-\stick\delta_i), & \delta_i>0,\ s_{i+1}\in\mathcal V.
	\end{cases}
\]
For the intervals of $\rec$ we take $\mathcal V=\emptyset$. In the last two cases, $h$ alone is the
density in time of the zero that ends the interval, as in Proposition~\ref{prop:zeros}, while $\stick h$ is the
weight of the atom at position $0$ at a prescribed time, as in the third case of Proposition~\ref{prop:kernel}. The
reading uses the skeleton alone, not the way its times arose. A recorded zero of $\rec_{\mathcal V}$ that is not a point
of $\rec$ is read as a zero. A point of $\rec$ at position $0$ is read as a zero for the intervals of $\rec_{\mathcal V}$
inside a constant interval of $\rec$, and for the interval of $\rec_{\mathcal V}$ that it ends inside an increasing interval
of $\rec$. All other ends are read as fixed. The left end of an increasing interval of $\rec$ is therefore fixed, since
there the middle bridge starts from $0$. A point of $\rec$ at position $0$ that ends an increasing interval is read as a
zero even when its time was a query time, at which the original process has an atom. This reading is the same in the
numerator and in the denominator of Proposition~\ref{prop:record}, and Lemma~\ref{lem:mixture} accounts for the factor
$\stick$ that it moves when one skeleton refines another. The reference measure $\nu_{\rec,\mathcal V}$ is Lebesgue
measure in the times of the new recorded zeros, and
$\big(\delta_0(\md z_v)+\md z_v\big)\big(\delta_0(\md\delta_v)+\md\delta_v\big)$ at each time $v\in \mathcal V$, where $z_v$ is
the position at $v$ and $\delta_v$ the local time gained since the previous time of $\mathcal T_\rec\cup \mathcal V$.

\begin{proposition}\label{prop:record}
Let $\rec\in\Recs$ or $\Recs_T$, and let $\mathcal V$ be a finite set of times that are not times of $\rec$, in $(0,T)$ for
$\rec\in\Recs_T$. If $\rec$ and $\rec'$ have $m$ and $m'$ intervals with respective factors $\varphi_i$ and $\varphi'_i$, then
\[
	\Ker_\rec\{\rec_{\mathcal V}\in\md\rec'\}~=~\frac{\prod_{i=0}^{m'-1}\varphi'_i}{\prod_{i=0}^{m-1}\varphi_i}\,\nu_{\rec,\mathcal V}(\md\rec').
\]
\end{proposition}

For $\rec=\{(0,x,0)\}$ the denominator is empty, and the proposition gives the density $\prod\varphi_i$ of $\rrec(\mathcal T)$
under $\Probx{x}$ for every finite $\mathcal T\ni0$. Integrating out the recorded zeros by Proposition~\ref{prop:zeros} leaves the
product of the transition densities $\pYl$ over consecutive times of $\mathcal T$, which is the finite-dimensional law of
$(\Y,\locT)$ on $\mathcal T$.

\begin{proof}[Proof of Proposition~\ref{prop:record}]
Consecutive times of $\mathcal T_\rec\cup \mathcal V$ lie in one interval of $\rec$, or both after its last point $s_m$. The points that
\eqref{eq:revealed} adds between them are functions of the path on that interval, or after $s_m$. At an end $s_i$ of an
interval of $\rec$ the skeleton adds no zero, because either the interval is constant or $z_i=0$. The points of $\rec_{\mathcal V}$ in
$[s_i,s_{i+1}]$, shifted by $(s_i,a_i)$, are therefore the skeleton of the shifted path at the shifted times, and the
same holds after $s_m$. The shifted paths are independent under $\Ker_\rec$, so the law of $\rec_{\mathcal V}$ is the product of
the laws of these skeletons. It remains to find the law of the skeleton of one shifted path. Let
$\mathcal V'=\{v_1<\dots<v_n\}$ be the shifted times in one interval of length $\Delta$, or after $s_m$.

Under $\Bes_{z,z',\Delta}$ the local time is $0$, so every pair is constant and the skeleton holds only the points at
$\{0\}\cup \mathcal V'\cup\{\Delta\}$. Their positions have the density of Lemma~\ref{lem:bridges-paths}(a), which is the product
of the factors of the constant intervals divided by the factor of $[0,\Delta]$. An end at position $0$ is an end at a zero,
and its factor is $h$ or $2h$, as the Bessel(3) bridge from or to $0$ starts or ends at a zero.

Under $\Mid_{\delta,t}$, with $\Delta=t+\stick\delta$, the local time at the end is $\delta$, and $\locT_v<\delta$ for
$v<\Delta$ (Lemma~\ref{lem:bridges-paths}(c)). The pair $(v_n,\Delta)$ is therefore increasing. It records the first zero
$d_{v_n}$ if $\Y_{v_n}>0$, and no last zero, since $\Y_\Delta=0$. Let $e\ass d_{v_n}$ if $\Y_{v_n}>0$ and $e\ass v_n$
otherwise. Then $\Y_e=0$ and $\locT_e=\locT_{v_n}$, and $e<\Delta$. For a path $w$, let $\Phi(w)$ be the set
$\rrec(\{0\}\cup\mathcal V')(w)\cup\{(e,0,\locT_e),\,(\Delta,0,\delta)\}$, the skeleton read from the path up to $e$ with the
prescribed endpoint appended. It is a function of the path stopped at $e$, and under $\Mid_{\delta,t}$ it equals
$\rrec(\{0\}\cup\mathcal V'\cup\{\Delta\})$. In the operational time of the first passage bridge, $e$ is the stopping time
$\varsigma\ass\sigma_e$. Lemma~\ref{lem:bridges-paths}(c) at $\varsigma$, with $\BM_\varsigma=\locT_e$, gives for every
bounded measurable $G$
\[
	\Mid_{\delta,t}\big[G\big(\rrec(\{0\}\cup \mathcal V'\cup\{\Delta\})\big)\big]
	~=~\frac1{h(\delta,t)}\,\E_0\Big[G(\Phi)\,\id_{\{\locT_e<\delta\}}\,
	h\big(\delta-\locT_e,\ \Delta-e-\stick(\delta-\locT_e)\big)\Big] .
\]
The last factor is the factor of the increasing interval $(e,\Delta)$, and $h(\delta,t)$ is the factor of $[0,\Delta]$.
Both factors read the right end $\Delta$ in the same way, so the choice between $h$ and $\stick h$ is the same in the
numerator and in the denominator, and it cancels. It remains to find the law of $\rrec(\{0\}\cup \mathcal V')$ under $\Probx{0}$, and of the
first zero after $v_n$, which is the case of the free path below. Given $\Y_{v_n}=z_n>0$, the time $d_{v_n}-v_n$ is the
first passage of a Brownian motion from $z_n$ to $0$, by Lemma~\ref{lem:markov} at $v_n$, and it has the density
$h(z_n,\iv)$ by \eqref{eq:h}. This is the factor of the constant interval $(v_n,d_{v_n})$.

Under $\Probx{z}$, the Markov property \eqref{eq:markov} at each time of $\{0\}\cup\mathcal V'$ (Lemma~\ref{lem:markov})
factors the joint density of the values and the recorded zeros over the intervals between consecutive times. On each
interval, after the shift of the local time by its value at the left end, the factor is the density of
Proposition~\ref{prop:kernel} or of Proposition~\ref{prop:zeros} from the left value. Let $u<u'$ be consecutive times.
If $\locT_{u'}=\locT_u$, the interval $[u,u']$ contributes the density $p_{\mathfrak A}(\Y_{u'},u'-u\mid\Y_u)$ of
Proposition~\ref{prop:kernel}, the factor of one constant interval with fixed ends. If $\locT_{u'}>\locT_u$, Proposition~\ref{prop:zeros} for the copy from $\Y_u$
over the time $u'-u$ gives the joint density of $(d_u,\,g_{u'},\,\Y_{u'},\,\locT_{u'}-\locT_u)$. When $\Y_u,\Y_{u'}>0$, it
is the product of the factors of the intervals $(u,d_u)$, $(d_u,g_{u'})$ and $(g_{u'},u')$. When $\Y_u=0$, the factor
$h(\Y_u,\iv)$ is absent and $d_u=u$. When $\Y_{u'}=0$, the factor $\stick$ replaces $2h(\Y_{u'},\iv)$, the increasing
interval ends at the fixed time $u'$, and $g_{u'}=u'$.

The factors of the intervals of $\rec_{\mathcal V}$ inside one interval of $\rec$, divided by the factor of that interval, are the
three laws just found, and the free path after $s_m$ contributes its factors without a denominator. This is the
display.
\end{proof}

\subsection{The bridge as a law of paths}\label{app:bridge-path}

The bridge of Section~\ref{sec:bridge} is a ratio of densities of finitely many
values, and Theorem~\ref{thm:refine} under the bridge needs it as a law of paths. The bridge draws its first and last
zeros, and the path given them is a Bessel(3) bridge, a middle bridge and a Bessel(3) bridge. Let $x,y\geq0$, $r\geq0$
and $T>0$ with $\pYl(y,r,T\mid x)>0$. For $r>0$, we write $\bar\rec(f,g)$ for the completed input of
Section~\ref{sec:repeated}, which consists of $(0,x,0)$, the first zero $(f,0,0)$ if $x>0$, the last zero $(g,0,r)$ if
$y>0$, and $(T,y,r)$, and we define
\begin{equation}\label{eq:bridge-path}
	\Probx{x}\{\iv\mid\Y_T=y,\ \locT_T=r\}~\ass~\int\Ker_{\bar\rec(f,g)}\;\Probx{x}\{\hT_0\in\md f,\ g_T\in\md g\mid\Y_T=y,\ \locT_T=r\},
\end{equation}
a law on $\mathcal C[0,T]$, where the law of $(f,g)$ is that of Corollary~\ref{cor:zeros}. For $r=0$, which requires
$x,y>0$, the bridge has no zero, and we define $\Probx{x}\{\iv\mid\Y_T=y,\ \locT_T=0\}\ass\Bes_{x,y,T}$, the kernel
of the skeleton $\bar\rec\ass\{(0,x,0),(T,y,0)\}$.

\begin{proposition}\label{prop:bridge-path}
Let $x,y\geq0$, $r\geq0$ and $T>0$ with $\pYl(y,r,T\mid x)>0$.
\begin{enumerate}
	\item[(a)] Let $\mathcal T$ be a finite set of times in $[0,T]$
containing $0$ and $T$. Under the bridge with endpoint $(y,r)$,
the skeleton $\rrec(\mathcal T)$ has density
\[
    \frac{\prod_i\varphi_i}{\pYl(y,r,T\mid x)}
\]
with respect to the following measure. The positions at the interior
times $v\in\mathcal T\setminus\{0,T\}$ have reference measure
$\delta_0(\md z_v)+\md z_v$, and each recorded zero whose time is
not prescribed contributes Lebesgue measure in its time.

For the local times, sum over the possible patterns of equalities
among the levels at the times of $\mathcal T$, with initial level
$0$ and terminal level $r$. For $r>0$, write the distinct levels
of a given pattern as
\[
    0=b_0<b_1<\cdots<b_k=r .
\]
The local-time measure on this pattern is
$\md b_1\cdots\md b_{k-1}$ on the indicated range, with mass one
when $k=1$. For $r=0$ there is a single pattern, with every level
zero and mass one. The local time at each recorded zero is supplied
by its neighbouring observation. In particular, the finite-dimensional
laws of the bridge are the ratios \eqref{eq:bridge-def} with
$\mathfrak w=(\Y_T,\locT_T)$.
	\item[(b)] For every bounded measurable $F$ on $\mathcal C[0,T]$ and every bounded measurable $\phi$,
	\[
		\begin{aligned}
		\E_x\big[F\big((\Y,\locT)|_{[0,T]}\big)\,\phi(\Y_T,\locT_T)\big]
		~=~\int&\E_x\big[F\bigm|\Y_T=y,\ \locT_T=r\big]\,\phi(y,r)\,\pYl(y,r,T\mid x)\\
		&\cdot\big(\delta_0(\md y)+\md y\big)\big(\delta_0(\md r)+\md r\big).
		\end{aligned}
	\]
	\item[(c)] For $x=y=0$ and $r>0$, \eqref{eq:bridge-path} is the middle bridge $\Mid_{r,\,T-\stick r}$.
\end{enumerate}
\end{proposition}

Part (b) says that \eqref{eq:bridge-path} is a regular version of the conditional law of the path given its endpoint.
Part (a) says that it agrees with the ratios \eqref{eq:bridge-def} at every endpoint, not only at almost every one.
Part (c) identifies the middle bridge, which Section~\ref{sec:bridge-given-zeros} defines by conditioning on $\hT_r=t$,
with the bridge from $(0,0)$ to $(0,r)$.

\begin{proof}
(a) If $r=0$, then $x,y>0$ and the bridge is $\Bes_{x,y,T}$.
Every prescribed level is zero. Lemma~\ref{lem:bridges-paths}(a)
gives the stated density in the interior positions, with denominator
$p_{\mathfrak A}(y,T\mid x)=\pYl(y,0,T\mid x)$.

Let $r>0$. By Proposition~\ref{prop:zeros} and
Corollary~\ref{cor:zeros}, the completed skeleton $\bar\rec(f,g)$
at $\{0,T\}$ has density $\prod_i\varphi_i/\pYl(y,r,T\mid x)$
in its free zero times. Given $\bar\rec(f,g)$, the path has law
$\Ker_{\bar\rec(f,g)}$. Its first and last zeros on $[0,T]$ are
$f$ and $g$ almost surely, since the outer Bessel(3) bridges are
positive in their interiors. Thus the final skeleton determines
$\bar\rec(f,g)$, and Lemma~\ref{lem:path-records}(c) gives the following equality almost surely.
\[
    \rrec(\mathcal T)
    ~=~\bar\rec(f,g)_{\mathcal T\setminus\{0,T\}}
    .
\]
Proposition~\ref{prop:record} gives the density of this refinement
as a ratio of interval factors. Multiplication by the density of
$\bar\rec(f,g)$ cancels its interval factors. When $y=0$, the
initial density has the factor $\stick h$ at the prescribed time
$T$. The refinement ratio uses $h$ at this inherited endpoint in
both numerator and denominator, so the product retains $\stick h$
for the final interval ending at $T$. The resulting density is
$\prod_i\varphi_i/\pYl(y,r,T\mid x)$.

It remains to identify the coordinates in which this is a density.
The local time is zero up to $f$ and equals $r$ from $g$ to $T$.
At every prescribed time strictly between $f$ and $g$ its level
lies in $(0,r)$, by Lemma~\ref{lem:bridges-paths}(c). Fix the
pattern of equalities among the prescribed levels. At the first
observation of each distinct interior level $b_j$, the gain from
the preceding level is $b_j-b_{j-1}>0$ and contributes Lebesgue
measure. Each repetition of that level contributes the unit atom
at zero gain. Hence the free gains and levels satisfy
\[
    \md(b_1-b_0)\cdots\md(b_{k-1}-b_{k-2})
    ~=~\md b_1\cdots\md b_{k-1}
    \qquad (k\geq2),
\]
by the triangular change of variables with unit diagonal.
The terminal level $r$ is already fixed, even when an interior
observation is the first prescribed time at that level. For $k=1$
there is no free level. The interior positions retain their
reference measures, and each free recorded zero time occurs once,
either in $\bar\rec(f,g)$ or in its refinement. This gives the
measure in (a). Integrating out the recorded zeros gives the
finite-dimensional ratio \eqref{eq:bridge-def}, by
Proposition~\ref{prop:zeros} on each interval between consecutive
prescribed times.

(b) First let $F$ depend only on $\rrec(\mathcal T)$. Under
$\Probx{x}$, Proposition~\ref{prop:record} gives the density
$\prod_i\varphi_i$ in the successive gains, positions and free
recorded zero times. Fix a local-time pattern with at least one
positive gain, and write $0=b_0<b_1<\cdots<b_k$ for its distinct
levels. In chronological order its positive gains are
$b_1-b_0,\ldots,b_k-b_{k-1}$. Changing from these gains to the
levels gives
\[
    \md b_1\cdots\md b_k
    ~=~(\md b_1\cdots\md b_{k-1})\,\md r,
    \qquad r=b_k .
\]
The factor in parentheses is precisely the local-time measure in
(a). On the all-zero pattern, the endpoint measure instead supplies
$\delta_0(\md r)$ and there is no free level. Summing the patterns
therefore recovers the unconditioned local-time measure. Integrating
the fixed terminal position against $\delta_0(\md y)+\md y$ likewise
restores its position measure.

Finally, the factor $\pYl(y,r,T\mid x)$ in (b) cancels the
denominator in (a). Both sides of (b) thus integrate $F\,\phi$
against the same skeleton law. In particular the identity holds
for bounded functions of finitely many path coordinates. These
coordinates generate the Borel $\sigma$-field on $\mathcal C[0,T]$,
so the monotone class theorem gives (b) for every bounded
measurable $F$.

(c) For $x=y=0$ the skeleton $\bar\rec$ has the two points $(0,0,0)$ and $(T,0,r)$, one increasing interval and no
recorded zeros besides its ends. By the definition of $\Ker$, the path is $\Mid_{r,T-\stick r}$.
\end{proof}

The kernel $\Ker$ is measurable in the endpoint, so \eqref{eq:bridge-path} is a kernel in $(y,r)$ and the integral in
(b) is defined. Under \eqref{eq:bridge-path} almost every path lies in $\Paths$, by Lemma~\ref{lem:kernel}.

\section{Exactness of refinement}\label{app:refine-proof}

We prove the two halves of the one-query identity (Lemmas~\ref{lem:query-law} and~\ref{lem:mixture}), which give
Proposition~\ref{prop:path-query}, and then Theorem~\ref{thm:refine} by induction on the number of queries.

\subsection{The law of the new skeleton}\label{app:query-law}

The query draws, at every skeleton and every time, the law of the skeleton of the path at the new time.

\begin{lemma}\label{lem:query-law}
Let $\rec\in\Recs$ or $\Recs_T$, and let $u$ be a time that is not a time of $\rec$, in $(0,T)$ for $\rec\in\Recs_T$.
Then
\[
	\Ker_\rec\{\rec_u\in\md\rec'\}~=~\qry_u(\rec,\md\rec').
\]
In addition, for $\qry_u(\rec,\iv)$-almost every $\rec'$, the skeleton $\rec'$ lies in $\Recs$ or $\Recs_T$, and
$\rec_u(w)=\rec'$ for $\Ker_{\rec'}$-almost every $w$.
\end{lemma}

\begin{proof}
By Proposition~\ref{prop:record} with $\mathcal V=\{u\}$, the density of $\rec_u$ is a ratio in which all factors cancel except
those of the interval of $\rec$ that contains $u$, or of the path after the last point. We compare this ratio with the
law that the query draws.

Let the interval $(s_i,s_{i+1})$ be constant. The path has no zero inside it, so the new skeleton splits it into two
constant intervals and adds no recorded zero. The ratio is
\begin{align*}
	&\frac{p_{\mathfrak A}(z,u-s_i\mid z_i)\,p_{\mathfrak A}(z_{i+1},s_{i+1}-u\mid z)}{p_{\mathfrak A}(z_{i+1},\Delta\mid z_i)}
	&&z_i z_{i+1}>0,\\
	&\frac{p_{\mathfrak A}(z,u-s_i\mid z_i)\,h(z,s_{i+1}-u)}{h(z_i,\Delta)}&&z_{i+1}=0,\\
	&\frac{2h(z,u-s_i)\,p_{\mathfrak A}(z_{i+1},s_{i+1}-u\mid z)}{2h(z_{i+1},\Delta)}&&z_i=0.
\end{align*}
These are the three laws of Lemma~\ref{lem:bessel}, which the query draws.

Let the interval be increasing. Its right end $s_{i+1}$ has the same factor in the numerator and in the denominator. If
$z=0$, the new skeleton splits the interval into two increasing intervals, and the first one ends at the fixed time $u$.
The ratio is
\[
	\frac{\stick h(a,\ u-s_i-\stick a)\,h\big(\delta-a,\ s_{i+1}-u-\stick(\delta-a)\big)}{h(\delta,\ \Delta-\stick\delta)},
\]
which is the first display of Proposition~\ref{prop:middle} with $(r,t,s)=(\delta,\ \Delta-\stick\delta,\ u-s_i)$. If
$z>0$, the new skeleton splits the interval into four intervals. These are an increasing interval that ends at the
recorded zero $g_u$, the two constant intervals $(g_u,u)$ and $(u,d_u)$, and an increasing interval that starts at the
recorded zero $d_u$. The ratio is
\[
	\frac{h(a,\ g_u-s_i-\stick a)\;2h(z,\ u-g_u)\;h(z,\ d_u-u)\;h\big(\delta-a,\ s_{i+1}-d_u-\stick(\delta-a)\big)}{h(\delta,\ \Delta-\stick\delta)} .
\]
By \eqref{eq:h-conv}, its integral over $g_u$ and $d_u$ is the second display of Proposition~\ref{prop:middle}. Given
$(z,a)$, the remaining factors in $u-g_u$ and in $d_u-u$, normalized, are the two passage times conditioned on later
ones \eqref{eq:ptb} that the query draws.

After the last point, the ratio is the density of the skeleton of $\Probx{z_m}$ at the times $\{0,u-s_m\}$. The new pair has
the transition density from $z_m$ over the time $u-s_m$, which Algorithm~\ref{alg:forward} draws (Theorem~\ref{thm:forward}). Given
the pair, the ratio of Proposition~\ref{prop:zeros} to Proposition~\ref{prop:kernel} is Corollary~\ref{cor:zeros}, which
the completion draws.

The query draws its new values from densities, so almost surely the new local time $a$ satisfies $0<a<\delta$ in an
increasing interval, and after the last point the new pair has $a>0$ or $z_m+z>0$. The new zeros are ends of increasing
intervals, and a new point at position $0$ lies between two increasing intervals. So $\rec'\in\Recs$, or $\Recs_T$.
Under $\Ker_{\rec'}$ the path passes through the points of $\rec'$. It is positive on $(g,u]$ and on $[u,d)$, where
$g<u<d$ are the new recorded zeros, because Bessel(3) bridges are positive except at an end at $0$. After the last point,
it is positive on $[s_m,d)$ and on $(g,u]$ for the zeros $d$ and $g$ of the completion. So the zeros of the path are the
new recorded zeros of $\rec'$. The path lies in $\Paths$ by Lemma~\ref{lem:kernel}, and $\rrec(\mathcal T_\rec)(w)=\rec$, so
Lemma~\ref{lem:path-records}(b) gives $\rec_u(w)=\rec'$.
\end{proof}

\subsection{The path given the new skeleton}\label{app:mixture}

Lemma E.2 is the second half of the one-query identity. A query at \(u\) reveals a finer skeleton. 
Averaging the path law assigned to that skeleton over the possible answers to the query recovers 
the path law before the query.

\begin{lemma}\label{lem:mixture}
Let $\rec\in\Recs$ and $u>0$ in the unconditioned case, or let $\rec\in\Recs_T$ and
$0<u<T$ in the bridge case. Suppose that $u\notin\mathcal T_\rec$. Then
\[
    \Ker_\rec=\int \Ker_{\rec'}\,\qry_u(\rec,\md\rec').
\]
\end{lemma}

\begin{proof}
Both sides are probability laws on continuous paths by Lemmas~\ref{lem:kernel}
and~\ref{lem:query-law}. Fix a finite set $\mathcal W$ of additional times in the path
domain, disjoint from $\mathcal T_\rec\cup\{u\}$, and put $\mathcal V=\mathcal W\cup\{u\}$.
We compare the laws of the refined skeleton $\rec''=\rec_{\mathcal V}$.

Under the left side, Proposition~\ref{prop:record} gives the density
$D_\rec(\rec'')$ with respect to $\nu_{\rec,\mathcal V}$, where $D$ denotes the ratio of
interval factors in that proposition. Under the right side, first draw $\rec'$ by
querying at $u$, then read the skeleton at $\mathcal W$ from a path of law $\Ker_{\rec'}$.
By Lemma~\ref{lem:query-law} the path recovers $\rec'$ at $\mathcal T_\rec\cup\{u\}$,
and Lemma~\ref{lem:path-records}(c) gives $\rec_{\mathcal V}=\rec'_{\mathcal W}$
almost surely. Conversely, to recover $\rec'$ from $\rec''$, retain the observations at
$\mathcal T_\rec\cup\{u\}$. Between consecutive retained times with different levels,
read from $\rec''$ the first and last zeros required by \eqref{eq:revealed}.
Thus the two-stage calculation uses the density $D_\rec(\rec')D_{\rec'}(\rec'')$
against the iterated measures $\nu_{\rec,\{u\}}$ and $\nu_{\rec',\mathcal W}$.
New recorded zero times have densities, so coincidences with times of $\mathcal W$
may be discarded.

The interval factors of $\rec'$ cancel, with one endpoint convention to track.
If an increasing interval ends at $u$ with $\Y_u=0$, its first-stage factor is
$\stick h_{\rm old}$ because $u$ is a prescribed time. In the second stage $u$ is an
inherited endpoint. The factors for the incoming intervals in the denominator and
numerator are then $h_{\rm old}$ and $h_{\rm new}$, respectively, each evaluated at
the gain and operational duration of its own interval. Their contribution to the product is
\[
    (\stick h_{\rm old})\frac{h_{\rm new}}{h_{\rm old}}=\stick h_{\rm new},
\]
as required at the prescribed time $u$ in the direct calculation. Hence the product
of the two densities is $D_\rec(\rec'')$ on the refinements produced in two stages.

It remains to compare the reference measures. The atom at zero need not occur at the
same query time in the two calculations. For example, take $\mathcal W=\{v\}$ and
suppose $g_u<v<u$ in a positive excursion revealed by a query inside an increasing
interval of $\rec$. Write $a>0$ for the gain from its old left endpoint to $u$.
In coordinate order $(v,u)$, the direct calculation gives the gain vector $(a,0)$.
The two-stage calculation gives $(0,a)$.
In the second calculation $g_u$ is already recorded when $v$ is queried. Both
calculations integrate over the one level shared by $v$ and $u$.

To see this for any $\mathcal W$, fix the order and roles of the points of $\rec''$,
which prescribed positions are zero, and which local-time levels are equal,
including equalities with the levels of $\rec$. These choices partition the refinements
into finitely many cases. On each case, read the reference measure of
Proposition~\ref{prop:record} for a single refinement as follows. An observation in
a constant interval repeats the level at its left end. Inside an increasing interval
the level lies strictly between the endpoint levels almost surely, since the middle
bridge gains local time in every neighbourhood of both ends
(Lemma~\ref{lem:bridges-paths}(c)). Beyond the last point, the levels are
nondecreasing from its level. Thus any observation at an old level has an earlier
point of the starting skeleton at that level. Among observations at a new level,
only the first in time contributes a positive gain. Every repetition contributes the
unit mass at zero. Recorded zeros copy the level of an adjacent observation.

Consequently the free local-time coordinates of one refinement are precisely the
distinct observed levels absent from its starting skeleton. Indeed, in order of
introduction write $b_j$ for such a level and $r_j$ for its positive gain. Then
$b_j=c_j+r_j$, where the predecessor level $c_j$ is either fixed or has already
been introduced. In these free coordinates the change of variables is triangular
with unit diagonal. Thus, if there are $k$ new levels, the local-time factor on
their allowed range is
\[
    \md r_1\cdots\md r_k=\md b_1\cdots\md b_k .
\]
When there are no new levels, all gains at the query times are zero and the
local-time factor has mass one.

Apply this description first with starting skeleton $\rec$ and query set $\mathcal V$.
Now apply it with query set $\{u\}$ and then with starting skeleton $\rec'$ and query
set $\mathcal W$. The first stage introduces only $\locT_u$, if that level is new.
The second introduces exactly the levels still absent from $\rec'$. Thus each
calculation integrates once over every distinct new level of $\rec''$. The positions
use the same atom or Lebesgue coordinate at each prescribed time, and the times of
the new recorded zeros occur once each, in either order. Since $\rec'$ is recovered
from $\rec''$, there is no multiplicity on passing to these common coordinates.

The two-stage law therefore agrees with the direct density on the cases it produces.
It has total mass one, as does the direct law, so those cases exhaust the direct law
up to a null set. The values at $\mathcal V$ are read from $\rec''$, and both paths
pass through the prescribed points of $\rec$. Since $\mathcal W$ was arbitrary, the
finite-dimensional laws agree. The coordinate maps generate the Borel $\sigma$-field
on the continuous path space (Appendix~\ref{app:path-record}), which proves the identity.
\end{proof}

\subsection{The conditional path identity}\label{app:path-query}

\begin{proof}[Proof of Proposition~\ref{prop:path-query}]
If $u$ is the time of a point of $\rec$, then $\rec_u(w)=\rrec(\mathcal T_\rec)(w)=\rec$ for $\Ker_\rec$-almost every $w$, by
Lemma~\ref{lem:path-records}(a), and $\qry_u(\rec,\iv)$ is the point mass at $\rec$. Otherwise, by
Lemma~\ref{lem:mixture} and then Lemma~\ref{lem:query-law},
\begin{align*}
	\int F(w)\,G\big(\rec_u(w)\big)\,\Ker_\rec(\md w)
	~&=~\int\!\!\int F(w)\,G\big(\rec_u(w)\big)\,\Ker_{\rec'}(\md w)\,\qry_u(\rec,\md\rec')\\
	&=~\int\Ker_{\rec'}[F]\;G(\rec')\;\qry_u(\rec,\md\rec').
\end{align*}
\end{proof}

\subsection{Proof of the refinement theorem}\label{app:thm-refine}

The proof is an induction on the number of queries. At each step the path,
given everything sampled so far, has the law of the kernel of the current skeleton. The next query time is measurable with
respect to what has been sampled, and the disintegration theorem lets us evaluate Proposition~\ref{prop:path-query} at
this random time. We write $X$ for the path $(\Y,\locT)$ under $\Probx{x}$, or under the bridge \eqref{eq:bridge-path},
and $\E$ for the expectation under this law.

We first reduce the input. Under $\Probx{x}$ with a fixed set $\mathcal T_0=\{0,s_1,\dots,s_n\}$, $0<s_1<\dots<s_n$,
the completed input has the law of the skeleton after $n$ queries at $s_1,\dots,s_n$,
\[
	\rec_0~\law~\qry_{s_n}\big(\cdots\qry_{s_1}(\{(0,x,0)\})\cdots\big).
\]
Each of these queries is a forward step from the last point and a completion of the new pair. By Theorem~\ref{thm:forward} and the Markov property, the
values form the Markov chain of the transition densities $\pYl$, and the zeros of each pair are drawn independently by
Corollary~\ref{cor:zeros}, as in the completed input. The completed input draws all these zeros after the values,
while the queries draw the zeros of each pair before the next forward step. The order does not change the joint law,
since the zeros of a pair depend only on the values at its ends and are conditionally independent of the other pairs
given the values. On the side of the path, $\rrec(\mathcal T_0)$ is the skeleton after the
same queries from $\rrec(\{0\})=\{(0,x,0)\}$. It is therefore enough to prove the theorem for $\mathcal T_0=\{0\}$, with the
first $n$ functions $\chi_j$ constant. Under the bridge, $\mathcal T_0=\{0,T\}$ and $\rec_0=\bar\rec(f,g)$.

For every bounded measurable function $F$ of the path, we prove the following identities by induction on $k$.
\begin{align}
	\E\big[F(X)\bigm|\mathcal G_k\big]~&=~\Ker_{\rrec(\mathcal T_k)}[F],\label{eq:ind-path}\\
	\big(\rec_0,\dots,\rec_k,\ \mathfrak v_1,\dots,\mathfrak v_k\big)~&\law~\big(\rrec(\mathcal T_0),\dots,\rrec(\mathcal T_k),\ \mathfrak v_1,\dots,\mathfrak v_k\big).\label{eq:ind-law}
\end{align}
For $k=0$ under $\Probx{x}$, $\rrec(\mathcal T_0)=\{(0,x,0)\}$ is constant, $\Ker_{\{(0,x,0)\}}=\Probx{x}$, and
$\rec_0=\{(0,x,0)\}$. Under the bridge, $\rrec(\{0,T\})=\bar\rec(f,g)$ almost surely by the proof of
Proposition~\ref{prop:bridge-path}(a), and \eqref{eq:ind-path} is the definition \eqref{eq:bridge-path}. The input
$\rec_0$ is drawn from the same law of $(f,g)$.

Let $k\geq1$, and suppose that \eqref{eq:ind-path} and \eqref{eq:ind-law} hold for $k-1$. Put
$\mathcal G\ass\mathcal G_{k-1}\vee\sigma(\mathfrak v_k)$, $\rec\ass\rrec(\mathcal T_{k-1})$ and $u\ass u'_k=\chi_k(\rec,\mathfrak v_k)$. Both
$\rec$ and $u$ are $\mathcal G$-measurable. The variable $\mathfrak v_k$ is independent of $X$ and of $\mathcal G_{k-1}$, so
\eqref{eq:ind-path} for $k-1$ still holds with $\mathcal G$ in place of $\mathcal G_{k-1}$, and $\Ker_\rec$ is a regular
conditional law of $X$ given $\mathcal G$. By Lemmas~\ref{lem:paths-full} and~\ref{lem:kernel}, $X\in\Paths$ almost
surely. The recorded zeros of $\rec$ are points of $\rrec(\mathcal T_{k-1}\cup\{u\})(X)$, by Lemma~\ref{lem:path-records}(c), and
the same lemma gives
\[
	\rrec(\mathcal T_k)~=~\rrec\big(\mathcal T_{k-1}\cup\{u\}\big)(X)~=~\rrec\big(\mathcal T_\rec\cup\{u\}\big)(X)~=~\rec_u(X).
\]

Let $F$ be a bounded measurable function of the path and $H$ a bounded measurable function of $k+1$ skeletons and $k$
numbers. We write $H(\rec')$ for $H\big(\rrec(\mathcal T_0),\dots,\rrec(\mathcal T_{k-1}),\rec',\mathfrak v_1,\dots,\mathfrak v_k\big)$, a function of
$\rec'$ whose other arguments are $\mathcal G$-measurable. By the disintegration theorem
\citep[Thm.~8.5(ii)]{kallenbergFoundationsModernProbability2021}, applied to the $\mathcal G$-measurable variables and the
regular conditional law $\Ker_\rec$, and since $(\rec,u,w)\mapsto\rec_u(w)$ is measurable (Lemma~\ref{lem:kernel}),
\begin{align*}
	\E\big[F(X)\,H\big(\rrec(\mathcal T_k)\big)\bigm|\mathcal G\big]~&=~\int F(w)\,H\big(\rec_u(w)\big)\,\Ker_\rec(\md w)\\
	&=~\int\Ker_{\rec'}[F]\;H(\rec')\;\qry_u(\rec,\md\rec').
\end{align*}
The second line is Proposition~\ref{prop:path-query}, which holds at every skeleton $\rec$ and every time $u$. With $F=1$
it says that the conditional law of $\rrec(\mathcal T_k)$ given $\mathcal G$ is $\qry_u(\rec,\iv)$. Applied to the function
$\rec'\mapsto\Ker_{\rec'}[F]\,H(\rec')$, this gives
\[
	\E\big[\Ker_{\rrec(\mathcal T_k)}[F]\,H\big(\rrec(\mathcal T_k)\big)\bigm|\mathcal G\big]~=~\int\Ker_{\rec'}[F]\;H(\rec')\;\qry_u(\rec,\md\rec').
\]
The two right sides agree. Taking expectations, $\E[F(X)\,H]=\E[\Ker_{\rrec(\mathcal T_k)}[F]\,H]$ for all such $H$, and every
bounded $\mathcal G_k$-measurable variable has the form $H\big(\rrec(\mathcal T_k)\big)$. This is \eqref{eq:ind-path} for $k$.

For \eqref{eq:ind-law}, the variable $\mathfrak v_k$ is independent of the past on both sides, so \eqref{eq:ind-law} for $k-1$
extends to $(\rec_0,\dots,\rec_{k-1},\mathfrak v_1,\dots,\mathfrak v_k)$. Both sides choose the next time by the same function
$\chi_k$. The algorithm draws $\rec_k=\qry_{u_k}(\rec_{k-1})$ with new elementary random variables, so its conditional
law given the past is $\qry_{u_k}(\rec_{k-1},\iv)$. On the side of the path, the conditional law of $\rrec(\mathcal T_k)$ given
$\mathcal G$ is $\qry_{u'_k}(\rrec(\mathcal T_{k-1}),\iv)$, as just shown. The two joint laws are therefore equal, and this is
\eqref{eq:ind-law} for $k$.

The first display of Theorem~\ref{thm:refine} is \eqref{eq:ind-law}, and the second follows by reading off the points at
the query times. The last statement of the theorem is \eqref{eq:ind-path}.

\bibliographystyle{sticky}
\bibliography{zotero,bibl,additions}
\end{document}